\documentstyle{amsart}

\newtheorem{thm}{Theorem}[section]
\newtheorem{prop}[thm]{Proposition}
\newtheorem{lem}[thm]{Lemma}
\newtheorem{cor}[thm]{Corollary}

\theoremstyle{definition}
\newtheorem{defn}[thm]{Definition}
\newtheorem{remark}[thm]{Remark}
\newtheorem{example}[thm]{Example}

\numberwithin{equation}{section}

\newcommand{\codim}{\operatorname{codim}}

\newcommand{\UD}{\operatorname{UD}}

\newcommand{\Spec}{\operatorname{Spec}}

\newcommand{\bbA}{{\Bbb A}}

\newcommand{\bbZ}{{\Bbb Z}}

\newcommand{\bbQ}{{\Bbb Q}}

\newcommand{\bbP}{{\Bbb P}}

\newcommand{\cO}{{\cal O}}

\newcommand{\fm}{{\frak m}}
\newcommand{\fp}{{\frak p}}

\newcommand{\isomo}{\overset{\sim}{=}}

\newcommand{\brokrarr}{\dasharrow}

\newcommand{\lf}{\mathopen}

\let\r=\mathclose

\newcommand{\notdiv}{\mathrel{\not|}}

\newcommand{\Ker}{\operatorname{Ker}}
\newcommand{\Aut}{\operatorname{Aut}}

\newcommand{\GL}{{\operatorname{GL}}}
\newcommand{\SL}{{\operatorname{SL}}}

\newcommand{\Sympl}{{\operatorname{Sp}}}
\newcommand{\PGL}{{\operatorname{PGL}}}
\newcommand{\PGLn}{{\operatorname{PGL}_n}}

\newcommand{\rank}{\operatorname{rank}}
\newcommand{\Stab}{\operatorname{Stab}}

\newcommand{\Char}{\operatorname{char}} 

\newcommand{\Gal}{\operatorname{Gal}}
\newcommand{\Galois}{\Gal}
\newcommand{\lra}{\longrightarrow}

\newcommand{\Mat}{{\operatorname{M}}}
\newcommand{\Mn}{\Mat_n}

\newcommand{\trdeg}{\operatorname{trdeg}}

\newcommand{\Oct}{\bold O} 

\let\tilde=\widetilde
\let\hat=\widehat

\let\to=\longrightarrow

\begin{document}
\title[Splitting fields]{Splitting fields of G-varieties}
\author[Z. REICHSTEIN and B. YOUSSIN,  8-30-99]
{Z. Reichstein and B. Youssin} 
\address{Department of Mathematics, Oregon State University,
Corvallis, OR 97331}
\thanks{Z. Reichstein was partially supported by NSF grant DMS-9801675}
\email{zinovy@@math.orst.edu}
\address{Department of Mathematics and Computer Science,
University of the Negev, Be'er Sheva', Israel\hfill\break
\hbox{{\rm\it\hskip\parindent Current mailing address\/}}: 
Hashofar 26/3, Ma'ale Adumim, Israel}
\email{youssin@@math.bgu.ac.il}
\subjclass{Primary: 14L30, 16K20, 20G10, 14E05.  Secondary: 14E15}

\begin{abstract}
Let $G$ be an algebraic group, $X$ a generically free $G$-variety, and
$K=k(X)^G$.
A field extension $L$ of $K$ is called a splitting field of $X$ if the image
of the class of $X$ under the natural 
map $H^1(K, G) \mapsto H^1(L, G)$ is trivial.
If $L/K$ is a (finite) Galois extension then $\Gal(L/K)$ is called
a splitting group of $X$. 

We prove a lower bound on the size of a splitting field of $X$ in terms
of fixed points of nontoral abelian subgroups of $G$. A similar result 
holds for splitting groups. We give a number of applications, including 
a new construction of noncrossed product division algebras.
\end{abstract}

\maketitle
\tableofcontents

\section{Introduction}

Let $k$ be an algebraically closed field of characteristic zero,
let $K$ be a finitely generated field extension of $k$ and let
$G$ is an algebraic group defined over $k$. 
Recall that elements of the non-abelian cohomology set $H^1(K, G)$ can be
identified with (birational classes of)
generically free $G$-varieties $X$ such that $k(X)^G = K$
(see \cite[Section 1.3]{popov}).
The set $H^1(K, G)$ has no group structure in general; 
however, $H^1(K, G)$ is equipped with a marked element, which we shall 
denote by $1$. This element is represented by the ``split" $G$-variety
$X \simeq X_0 \times G$, where $k(X_0) = K$ and
$G$ acts by left multiplication on the second factor.  A field extension
$L/K$ is said to be a splitting field for $u \in H^1(K, G)$ if 
$u \mapsto 1$ under the natural map $H^1(K, G) \lra H^1(L, G)$.

The nonabelian cohomology set $H^1(K, G)$ often allows 
a different interpretation:
its elements can be identified with certain algebraic objects 
defined over $K$, e.g, quadratic forms if $G = O_n$,
central simple algebras if $G = \PGL_n$, Cayley algebras, if $G = G_2$,
etc. These objects may be viewed as ``twisted forms" of a single ``split"
object. In such cases the above notion of a splitting field coincides
with the usual one.  We will review 
this interpretation of $H^1(K, G)$ in Section~\ref{sect.structured-spaces};
see also~\cite[Chapter~III]{serregc},~\cite[Chapter~X]{serrelf},
\cite[Section~29]{boi} or~\cite[Sections~6-8]{re}. 

Recall that a subgroup of $G$ is called {\em toral} if it lies in a torus 
in $G$.  Our main results on splitting fields are
Theorems~\ref{thm1.2a} and~\ref{thm1.2b}.

\begin{thm} \label{thm1.2a}
Let $X$ be a generically free primitive $G$-variety, $K = k(X)^G$, and
let $L/K$ be a splitting field for $X$.
Suppose $X$ has a smooth point fixed by a finite abelian $p$-subgroup 
$H$ of $G$. Then $[L:K]$ is divisible by $[H:H_T]$ for some 
toral subgroup $H_T$ of $H$.
\end{thm}

If $X$ is a generically free primitive $G$-variety, $K = k(X)^G$, 
and $L$ is a splitting field which is a (finite) Galois extension of $K$, 
then we shall refer to $\Gal(L/K)$ as a {\em splitting group} for $X$.  

\begin{thm} \label{thm1.2b}
Let $X$ be a generically free primitive $G$-variety and 
let $A$ be a splitting group for $X$. Suppose $X$ has 
a smooth point fixed by a finite abelian 
subgroup $H$ of $G$. Then $A$ contains an isomorphic
copy of $H/H_T$ for some toral subgroup $H_T$ of $H$.
\end{thm}

Our proofs of Theorems~\ref{thm1.2a} and~\ref{thm1.2b} are based on the
following results of~\cite{ry}: for any finite abelian subgroup $H$ of $G$,
the existence of a smooth $H$-fixed point on a (complete smooth)
$G$-variety $X$ is a birational invariant of $X$. Moreover, 
such points survive under dominant rational $G$-equivariant maps and under 
certain $G$-equivariant covers; see~\cite[Section 5 and Appendix]{ry}.
We review and further extend these results in Section~\ref{sect2}; see
Proposition~\ref{prop.proj} and Theorems~\ref{thm.resol},~\ref{going-down} 
and~\ref{going-up}.

Informally speaking, Theorem~\ref{thm1.2a} (respectively,
Theorem~\ref{thm1.2b}) may be viewed as a ``lower bound" on a
splitting field (respectively, a splitting group) of $X$.
In particular, if
$X$ is a vector space and $G$ acts linearly on $X$ then
$X$ has a smooth $G$-fixed point (namely, the origin) and, hence, in this 
case Theorem~\ref{thm1.2a} (respectively, Theorem~\ref{thm1.2b})
 can be applied to every 
finite abelian $p$-subgroup (respectively, subgroup) $H$ of $G$. Of course, 
Theorems~\ref{thm1.2a} and~\ref{thm1.2b} are only of interest
if $H$ is nontoral, since otherwise $H/H_T$ may be trivial.

Nontoral elementary abelian subgroups of algebraic groups have been 
extensively studied (see~\cite{bs}, \cite{borel}, \cite{steinberg}, 
\cite{serresubgr}); a complete classification was obtained
by Griess~\cite{griess}. To the best of our knowledge, 
nonelementary abelian subgroups have not been classified.
In Section~\ref{sect.examples} we apply Theorem~\ref{thm1.2a} 
to a number of specific groups $G$, where we have sufficient 
information about the depth of certain nontoral subgroups (see
Definition~\ref{def.depth}).
In particular, for $G = E_8$ we give a new proof 
of a theorem of Serre; see Corollary~\ref{cor.e-8}. 
Note, however, that the examples we give in Section~\ref{sect.examples}
are somewhat fragmentary, because we do not know any general results about
the depth of finite abelian subgroups in exceptional groups.
(Propositions~\ref{prop.ec^8} and~\ref{prop.e-7} represent our best efforts
in this direction; see also Corollary~\ref{cor:upper.depth}.) 
We hope that this question will attract the attention of group theorists
in the future, and that a more complete picture will emerge.

``Upper bounds" on the degrees of splitting fields, 
i.e., results of the form ``every $G$-variety can 
be split by a field extension of degree dividing $n(G)$", 
can be found in the paper~\cite{tits} of Tits. 
For a discussion of these results, including a table of values for $n(G)$, 
see Remark~\ref{rem.upper}.

In Sections~\ref{sect10} and \ref{sect.trdeg6} we apply Theorem~\ref{thm1.2b}, 
with $G = \PGLn$ to the theory of central simple algebras.
Recall that an element $\alpha \in H^1(K, \PGLn)$ may be (functorially)
identified with an $n^2$-dimensional central simple $K$-algebra $D_{\alpha}$; 
see Example~\ref{ex.pgl_n}. In particular, $L/K$ is a splitting field 
for $\alpha$ if and only if $L$ is a splitting field for $D_{\alpha}$, i.e.,
$D_{\alpha} \otimes_K L  = M_n(L)$. Recall that $D$ is an $H$-crossed product 
iff $H$ is a splitting group for $D$ and $|H| = \deg(D)$.

Let $\UD(n, k)$ be the universal division algebra of degree $n$, i.e.,
the division algebra generated by two generic matrices,
$X = (x_{ij})$ and $Y = (y_{ij})$, in $\Mn(k(x_{ij}, y_{ij}))$. 
Here $x_{ij}, y_{ij}$ are algebraically independent variables over $k$.
If the reference to $k$ is clear from the context, 
we shall write $\UD(n)$ in place of $\UD(n, k)$. A famous theorem of
Amitsur asserts that $\UD(n)$ is not a crossed product if $n$ is divisible
by $p^3$ for some prime $p$.

As an application of Theorem~\ref{thm1.2b} we will prove the following
result.

\begin{thm} \label{thm1.3} 
Let $Z(p^r)$ be the center of the universal division algebra $\UD(p^r)$,
let $K$ be a field extension of $Z(p^r)$ and let
$D = \UD(p^r) \otimes_{Z(p^r)} \, K$. Suppose $p^e$ is the highest
power of $p$ dividing $[K:Z(p^r)]$, where $e$ is a non-negative integer and
$e \leq r-1$. If $A$ is a splitting group for $D$ then 
\[ p^{2r - 2e -2} \mathrel| |A|  \, . \]
In particular, if $r \geq 2e + 3$  then $D$ is a noncrossed product.
\end{thm} 

If $K = Z(p^r)$, i.e., $D = \UD(p^r)$, we recover a theorem of 
Amitsur and Tignol; see~\cite[Theorem 7.3]{ta1}.
If $e = 0$, i.e., $D$ is a prime-to-$p$ extension of $\UD(p^r)$, we 
recover a theorem of Rowen and Saltman~\cite[Theorem~2.1]{rs} 
to the effect that $D$ is not a crossed product for any $r \geq 3$.

Abelian subgroups of $\PGLn$ carry a natural skew-symmetric form and
their nontoral subgroups are isotropic with respect to this form; see
Section~\ref{sect9}.
Thus symplectic modules and their Lagrangian submodules, used by Tignol 
and Amitsur to prove \cite[Theorem 7.3]{ta1}, naturally arise in our setting;
in particular, they will be used in the proof of Theorem~\ref{thm1.3} in
Section~\ref{sect10}.

It is likely that Theorem~\ref{thm1.3} can also be proved by 
an application of Amitsur's specialization technique, along 
the lines of~\cite[Section 2]{rs} and that such a proof will
go through in prime characteristic (assuming $p \notdiv \Char(k)$). 
We believe that our approach, based on the fixed points of nontoral
subgroups, is of independent interest; in particular, it shows that
Theorem~\ref{thm1.3} remains true if $\UD(p^r)$ is replaced by
any central simple algebra whose corresponding $\PGL_n$-variety 
has points fixed by certain nontoral subgroups of $\PGL_n$;
see Remark~\ref{rem1.3}.

\smallskip
As another application of Theorem~\ref{thm1.2b} with $G = \PGLn$, 
we construct a noncrossed product division algebra 
over a ``small" function field.  Since the time 
of Amitsur's original examples, two other noncrossed product constructions
have appeared in the literature, due, respectively, to 
Jacob---Wadsworth~\cite{jw} and Brussel~\cite{brussel}. Both of these 
examples have the property that their centers are ``smaller" and easier
to describe than the center of Amitsur's ``generic" example, $\UD(p^r, k)$.

The problem we address here is one of constructing noncrossed product
examples over ``small'' fields in the geometric setting, i.e.,
noncrossed products $D$ with center $K$ such that 
$K$ is a function field over an algebraically closed base
field $k$ of characteristic 0.  Moreover, we would like ``the size 
of $K$", as measured by $\trdeg_k(K)$, to be as small as possible. 

Note that $\trdeg_k(K)$ cannot be $\leq 1$ by Tsen's theorem.
Moreover, division algebras $D$ with $\trdeg_k (K) = 2$ 
are conjectured to be cyclic. 
At the other extreme, if $D = \UD(n)$ is Amitsur's 
original noncrossed product example (with $n$ divisible 
by $p^3$ for some prime $p$) then $\trdeg_k(K) = n^2 + 1$. 

In this paper we prove the following theorem.

\begin{thm} \label{thm.trdeg6} Let $p$ be a prime, $r \geq 2$ be an integer,
and $k$ be an algebraically closed field of characteristic 0. Then 
there exists a division algebra $D$ of degree $p^r$ with center $K$, 
such that

\smallskip
(a) $K$ is a finitely generated extension of $k$ of transcendence 
degree $6$ and

\smallskip
(b) no prime-to-$p$ extension of $D$ is a crossed product.
\end{thm}

The idea of the proof is as follows.
We show (see Section~\ref{sect10}) that it is enough to construct a
smooth $\PGL_{p^r}$-variety $X$ with two points
whose stabilizers are ``incompatible" symplectic modules
$(\bbZ/p^r\bbZ)^2$ and $(\bbZ/p\bbZ)^6$; such varieties are fairly
easy to construct.  The difficult 
part is to reduce the dimension of $X/\PGL_{p^r}$ to 6; this is done in 
Section~\ref{sect.trdeg6}.
Our argument there is based on a resolution result for
the fixed point loci of finite abelian subgroups 
(we show that the fixed-point set of a finite abelian subgroup $H$ can
be resolved in such way that it has a component of the minimal
possible codimension, equal to $\rank H$; see Theorem~\ref{thm.jun10a}) 
and on a form of Bertini's theorem in the equivariant setting 
(Theorem~\ref{thm.bertini}). We believe Theorems~\ref{thm.jun10a} 
and~\ref{thm.bertini} are of independent interest.

A. R. Wadsworth has pointed out to us 
that Theorem~\ref{thm.trdeg6} can be proved by modifying 
the arguments of~\cite{jw}. This approach, based on valuation theory, 
cohomology, and the Merkurjev---Suslin theorem, yields the desired 
result under the assumption that $p \notdiv \Char(k)$.

Throughout this paper we shall work over a fixed base field $k$ which 
will be assumed to be algebraically closed and  of characteristic 
zero. The assumption that $k$ should be algebraically closed is usually 
not essential: generally speaking, the problems we wish 
to consider (such as constructing noncrossed products or proving 
lower estimates on the size of splitting fields) can only become 
harder after passing to the algebraic closure. 
The characteristic zero assumption is more serious, since most 
of our proofs ultimately rely on canonical resolution of singularities
(via Proposition~\ref{prop.proj} and Theorem~\ref{thm.resol}).

\section*{Acknowledgements}
We would like to thank D. J. Saltman, J.-P. Serre, and A. R. Wadsworth
for helpful discussions.

\section{$G$-varieties}
\label{sect:G-var}

\subsection*{Preliminaries}
\label{sect2}

A $G$-variety $X$ is an algebraic variety with a $G$-action. 
Here $G$, $X$ and all other algebraic objects in this paper are assumed
to be defined over a fixed base field $k$. Unless otherwise specified,
we shall assume that $k$ is algebraically closed and of characteristic 0.
The $G$-action on $X$ is given by a morphism $G \times X \lra X$. 
If the reference to the action is clear from the context,
we shall write $gx$ for the image of $(g, x)$ under this map. Given $x \in X$,
the stabilizer of $x$ is defined as $\{g \in G \mid gx = x \}$; we will
denote this subgroup of $G$ by $\Stab_G(x)$ or simply $\Stab(x)$ if 
the reference to $G$ is clear from the context. 

By a morphism $X \lra Y$ of $G$-varieties, we shall mean a $G$-equivariant
morphism from $X$ to $Y$. The same goes for rational morphism, isomorphism,
birational morphism, etc., of $G$-varieties.

A $G$-variety $X$ is called {\em primitive} if $G$ transitively permutes 
the irreducible components if $X$. Equivalently, $X$ is primitive iff
$k(X)^G$ is a field.  Note that an irreducible $G$-variety is 
necessarily primitive and that the converse holds if $G$ is a connected group. 

If $X$ is a $G$-variety then any variety $Y$ with $k(Y) = k(X)^G$ is 
called a {\rm rational quotient variety} for $X$; we will often write
$Y = X/G$. Note that $X/G$ is only defined up to birational isomorphism and
that $X$ is a primitive $G$-variety iff $X/G$ is irreducible.
The inclusion $k(Y) = k(X)^G \hookrightarrow k(X)$ induces the rational 
quotient map $\pi \colon X \brokrarr X/G$. By a theorem of Rosenlicht 
$\pi^{-1}(y)$ is a single $G$-orbit for $y$ in general position in $X/G$
(see~\cite[Theorem 2]{rosenlicht1},~\cite{rosenlicht2}).

A $G$-variety $X$ is called {\em generically free} if $\Stab(x) = \{ 1 \}$
for $x$ in general position in $X$.  We will usually consider 
$G$-varieties that are both primitive and generically free.
Up to birational isomorphism,  a primitive generically free $G$-variety
may be viewed as a principal $G$-bundle over $X/G$ and thus represents 
a class in $H^1(K, G)$, where $K = k(X)^G$; see \cite[Section 1.3]{popov}.
We shall return to this connection in Section~\ref{sect.structured-spaces}.

It is often convenient to have a concrete (biregular) model for $X/G$. 
If $G$ is a finite group then, under rather mild assumptions on $X$, 
we have such a model in the form of a {\em geometric quotient}, which we shall
denote by $X//G$. Here the quotient map $\pi \colon X \lra X//G$ 
is regular and each fiber of this map is a single $G$-orbit. 
For a precise definition and a detailed discussion of the geometric 
quotient we refer the reader to \cite[Section 4.2]{pv}.

\begin{lem} \label{lem.quasi-proj}
Let $G$ be a finite group and $X$ be a normal quasiprojective $G$-variety.
Then

\smallskip
(a) $X$ is covered by affine open $G$-invariant subsets.

\smallskip
(b) There exists a geometric quotient map $\pi \colon X \lra X//G$. 

\smallskip
(c) Moreover, if $X$ is projective then so is $X//G$.
\end{lem}

\begin{pf} (a) By a theorem of Kambayashi, we may assume without loss 
of generality
that $X \subset \bbP(V)$, where $V$ is a finite-dimensional vector space,
and $G$-acts linearly on $X$, via a representation $G \lra \GL(V)$; 
see~\cite[Theorem 2.5]{kambayashi} or~\cite[Theorem 1.7]{pv}. 

We want to show that every 
$x \in X$ has an affine $G$-invariant neighborhood in $X$.
To construct this neighborhood, choose a homogeneous polynomial 
$h \in k[V]$ such that 
$h(gx) \neq 0$ for every $g \in G$ but $h(y) = 0$ for every 
$y \in \overline{X} - X$. After replacing $h$ by the product of $g^*h$ over
all $g \in G$, we may assume $h$ is $G$-invariant. Now $
\{ z \in \overline{X} \, | \, h(z) \neq  0 \}$
is a desired affine $G$-invariant neighborhood of $x$.

\smallskip
(b) Follows from part (a) and~\cite[Theorem 4.14]{pv}. 

\smallskip
(c) See~\cite[Theorem 4.16]{pv}.
\end{pf}

\subsection*{The variety $X_L$}

Let $X$ be a generically free primitive $G$-variety, let $K = k(X)^G$ and
let $cl(X)$ be the class of $X$ in $H^1(K, G)$. Suppose $L$ is a finitely 
generated field extension of $K$. Then $X_L$ is defined as the $G$-variety
representing the image of $cl(X)$ under the natural map 
$H^1(K, G) \lra H^1(L, G)$.  In other words, $cl(X) \mapsto cl(X_L)$ 
under this map.

To construct $X_L$ explicitly, let $Y \brokrarr X/G$ be a rational map
such that $k(Y)/k(X/G)$ is precisely the extension $L/K$. Note that such 
a rational map exists because $L$ is finitely generated over $K$ and, hence,
over $k$. Now we set $X_L = Y \times_{X/G} X$, where the $G$-action on $X_L$
is induced from the $G$-action on $X$; cf.~\cite[Section 2.6]{re}.

We emphasize that $X_L$ is only defined up to birational isomorphism
(of $G$-varieties). We we will often want to work with a specific model 
for $X_L$ which is smooth or projective or has ``small" stabilizers (or
all of the above).  The existence of such models is guaranteed 
by Proposition~\ref{prop.proj} and Theorem~\ref{thm.resol}.

\subsection*{Smooth projective models for $G$-varieties}

\begin{prop} \label{prop.proj} Every $G$-variety is birationally
isomorphic to a smooth projective $G$-variety.
\end{prop}

\begin{pf} Let $X$ be a $G$-variety. By~\cite[Proposition 7.1]{ry},
$X$ is birationally isomorphic to a complete $G$-variety. (Note that 
the proof of~\cite[Proposition 7.1]{ry} is based on Sumihiro's 
equivariant completion theorem.) Thus we may assume without loss 
of generality that $X$ is complete. 

Now by~\cite[Theorem 2.5]{kambayashi}
there exists a projective representation $G \hookrightarrow \PGL(V)$ 
and a closed $G$-invariant subvariety $X'$ of $\bbP(V)$ such that $X$ 
and $X'$ are birationally isomorphic as $G$-varieties. After replacing
$X$ by $X'$, we may assume $X$ is projective.  Now apply the canonical
resolution of singularities theorem
(see either~\cite[Theorem~7.6.1]{Vil2} or~\cite[Theorem~13.2]{bm})
to $X$ to construct a smooth projective model.
\end{pf}

\begin{defn} \label{def.levi}
We shall call an algebraic group $H$ {\em Levi-commutative} if
$H$ is a semidirect product of a diagonalizable group $D$ and
a unipotent group $U$, where $U \triangleleft H$ is the unipotent 
radical of the identity component $H_0$ of $H$.
\end{defn}

We shall denote $U=R_u(H_0)$.

\begin{lem} \label{lem.levi}
Let $H$ be an algebraic group and let $H_0$ be the identity component of $H$.
The following conditions are equivalent.

\smallskip
(i) $H$ is Levi-commutative, 

\smallskip
(ii) $H/R_u(H_0)$ is commutative,

\smallskip
(iii) every reductive subgroup of $H$ is commutative, and

\smallskip
(iv) every linear representation of $H$ has 1-dimesional $H$-invariant 
subspace.
\end{lem}

\begin{pf} The equivalence of (i), (ii) and (iii) follows from the Levi
decomposition theorem (see~\cite[Section 6.4]{ov}). The equivalence of (i) and
(iv) is proved in~\cite[Lemma A.1]{ry}.
\end{pf}

\begin{thm} \label{thm.resol} 
Let $X$ be a $G$-variety. Then there exists a sequence of blowups
\[ X_n \lra X_{n-1} \lra \dots \lra X_0 = X \] 
such that $X_n$ is smooth
and $\Stab(x)$ is Levi-commutative for every $x \in X_n$.
\end{thm}

\begin{pf} See~\cite[Theorem 1.1]{ry}.
\end{pf}

\subsection*{Going up and going down}

\begin{thm}[Going down]\label{going-down}
Let $H$ be a Levi-commutative group (see Definition~\ref{def.levi})
and let $X \brokrarr Y$
be a rational map of $H$-varieties. Suppose $Y$ is complete and 
$X$ has a smooth point fixed by $H$. Then $Y$ has a smooth point fixed by 
$H$.
\end{thm}

\begin{pf} See~\cite[Propositions 5.3 and A.2]{ry}.
\end{pf}

\begin{thm}[Going up]\label{going-up}
Let $H$ be a finite abelian group of prime power order $p^n$ and let 
$f \colon X \brokrarr Y$ be a dominant rational map of $H$-varieties. 
Suppose

\smallskip
(i) $Y$ is irreducible,

\smallskip
(ii) $X$ is complete,

\smallskip
(ii) $\dim(X) = \dim(Y)$ and $f$ is dominant,

\smallskip
(iv) $H$ has a smooth fixed point in $Y$, and 

\smallskip
(v) $\deg(f)$ is not divisible by $p^{e+1}$.

\smallskip\noindent
Then there exists a subgroup 
$H'$ of $H$ such that $|H'| \geq p^{n-e}$ and

\smallskip
(a) $H'$ has a fixed point on $X$. 

\smallskip
(b) Moreover, if $X \brokrarr Z$ is a rational map 
of $H$-varieties then $H'$ has a fixed point on $Z$. 
\end{thm}

This theorem is a generalization of~\cite[Propositions 5.5 and A.4]{ry}, 
where $e$ is assumed to be $0$.  Our proof 
below is based on an argument of Koll\'ar and Szab\'o; 
cf.\ \cite[Proposition~A.4]{ry}.  Our applications 
will only use (a); however, part (b) is needed for the inductive argument.

A convenient way to visualize the setting of Theorem~\ref{going-up} is
by means of the diagram
\[  \begin{array}{ccc}
X & - \stackrel{h}{-} \; \rightarrow & Z \\  
| \\
\llap{$f$ }|\\
\downarrow \\
Y
\end{array} \]
If $e = 0$ the theorem allows us to lift an $H$-fixed point from $Y$ to $X$,
then transport it to $Z$. (A  similar but slightly weaker statement is
true if $e \geq 1$.) We will make use of such diagrams in the proof.

\begin{pf}
The proof is by induction on $d=\dim Y$. 
If $\dim Y=0$ then $Y$ is a point, $X$ is a set of $\deg(X/Y)$ points,
and the desired result follows from a simple counting argument.

To perform the induction step, we assume that the theorem holds 
whenever $\dim(Y) \leq d-1$.  Let $y \in Y$ be a smooth fixed point, 
$B_y(Y)$ be the blowup of $Y$ at $y$, and $Y' \subset B_yY$ be
the exceptional divisor. Note that $H$ acts linearly on
$Y' = \bbP^{\dim(Y) -1}$ and, hence, has a fixed point in $Y'$;
see Lemma~\ref{lem.levi}(iv).  This fixed point will
be smooth because every point of $Y'$ is smooth.

Let $X=\bigcup_i X_{i}$ be the decomposition of $X$ into irreducible
components.
It is enough to find the required fixed point in one
of the components $X_i$ which is mapped dominantly onto $Y$; thus, we
may assume that all $X_i$ are mapped dominantly onto $Y$.

It follows that each map $X_i\brokrarr B_yY$ is dominant; let
$\overline X_i$ be the normalization of $B_yY$ in the field 
of rational
functions on $X_{i}$, $\overline X$ be the disjoint union of all
$\overline X_i$ (in other words, $\overline X$ is the normalization of
$B_yY$ in the ring $k(X)=\bigoplus_ik(X_{i})$), and
$\overline f:\overline X\to B_yY$ be the natural morphism. Clearly, $H$
acts on $\overline X$ and $\overline f$ is $H$-equivariant.
Each $\overline X_i$, it is birationally isomorphic
to $X_{i}$. Together these birational isomorphisms yield
an $H$-equivariant birational isomorphism $\overline X\brokrarr X$.

Let $F_1, F_2, \ldots\subset \overline X$ be
the divisors lying over $Y'$. Note that even though $\overline X$ is not
necessarily complete, each $F_i$ is complete since it is mapped
finitely to a complete variety $Y'$.

The group $H$ acts on the set  $\{F_i\}$.
Let ${\cal F}_j$ denote the $H$-orbits in $\{F_i\}$. Choose a divisor
$F^*_j\in {\cal F}_j$ in each orbit.
By the ramification formula (see, e.g., \cite[Corollary~XII.6.2]{lang}),
\[ \deg (X/Y)=\sum_j 
|{\cal F}_j|\cdot\deg (F^*_j/E)\cdot e(\overline f, F^*_j) \, ,\]
where $e(\overline f,F^*_j)$ denotes
the ramification index of $\overline f$ at the generic point of $F^*_j$.
Since $\deg (X/Y)$ is not divisible by $p^{e+1}$,  
$|{\cal F}_j|\cdot\deg (F^*_j/E)$ is not divisible by $p^{e+1}$ for some $j$.
For this $j$,
set $X' = HF_j^* = \bigcup_{F_i \in {\cal F}_j} F_i$; this variety
is complete since each $F_i$ is complete.  Let
$f' \colon X' \to Y'$ be the restriction of
$\overline f \colon \overline X \to B_z(Y)$ to $X'$; the degree of
$f'$ is equal to
$|{\cal F}_j|\cdot\deg (F^*_j/E)$, and hence, is not divisible by $p^{e+1}$.
Let $h' \colon X'\brokrarr X$ be the restriction of
$\overline X\brokrarr X$ to $X'$. Note that $h'$ is well-defined, since
$\overline X$ is normal and $X'$ is a divisor in $\overline X$.

By our construction, $\dim(Y') = d-1$ and conditions (i)--(v) hold for 
the map $f' \colon X' \brokrarr Y'$.  Applying the induction assumption to 
the diagram
\[  \begin{array}{ccc}
X' & - \stackrel{h'}{-} \; \rightarrow & X \\  
| \\
\llap{$f'$ }|\\
\downarrow \\
Y'
\end{array} \]
we prove part (a). Applying the induction assumption to the diagram
\[  \begin{array}{ccccc}
X' & \stackrel{h'}{\brokrarr} & X  & \stackrel{h}{\brokrarr} & Z \\  
| \\
\llap{$f'$ }|\\
\downarrow \\
Y'
\end{array} \]
we prove part (b).
\end{pf}

\begin{remark} \label{rem.char} Theorems~\ref{going-down} 
and~\ref{going-up} are valid 
over an algebraically closed field of arbitrary characteristic;
the proofs given above are characteristic-free.
\end{remark}

\section{Two interpretations of $H^1$}
\label{sect.structured-spaces}

Let $A$ be a finite-dimensional algebra over $k$. We do not assume
that $A$ is commutative, associative or has an identity element. Let 
$K$ is a field extension of $k$. We shall say that a $K$-algebra
$B$ is of type $A$ if $B_{\overline{K}} \simeq A_{\overline{K}}$, where
$\overline{K}$ is the algebraic closure of $K$. (Here, as usual,
$B_{\overline{K}} \stackrel{\rm def}{=} B \otimes_K \overline{K}$ and
$A_{\overline{K}} \stackrel{\rm def}{=} k \otimes_k \overline{K}$.)

Let $G = \Aut_k(A)$. It is easy to see that $G$ is a closed subgroup
of $\GL_{\dim(A)}$; thus it is an algebraic group. We now have the 
following bijections. 

\[ \begin{array}{clr} H^1(K, G) & \leftarrow \rightarrow & 
\left\{ \begin{array}{l} \text{Birational isomorphism classes} \\
 \text{of primitive generically free} \\
\text{$G$-varieties $X$ with $k(X)^G = K$} \end{array} \right\} \\
  \uparrow & &  \\
  \downarrow & &  \\
             & &  \\
\left\{ \text{$K$-algebras of type $A$} \right \}\\
\end{array}\] 
The horizontal correspondence is described in \cite[Section 1.3]{popov};
the vertical one in~\cite[Section X.2]{serrelf} 
or~\cite[Proposition 29.1]{boi}.  Note that the above 
diagram is functorial in $K$ and that the correspondences in it 
are bijections of pointed sets:
the identity elements of $H^1(K, G)$ corresponds to the ``split"
algebra $A_K$ and to the ``split" variety $X = X_0 \times G$, where 
$k(X_0)^G = K$; see Definition~\ref{def1.1}. 

In this paper we shall be primarily interested in passing 
back and forth between $G$-varieties and algebras of type $A$. 
In other words, we would like to construct an explicit correspondence
$X \mapsto B$ which completes the triangle in the above diagram. 
Given a generically free primitive $G$-variety $X$ with $k(X)^G = K$, 
we define $B = RMaps_{G}(X, A)$. Here $RMaps_{G}(X, A)$ is
the set of $G$-equivariant rational maps $X \brokrarr A$, where 
we view $A$ as a $k$-vector space with a $G$-action.
The $k$-algebra structure on $A$ gives rise to a $K$-algebra structure on
$B$ via
\[ \begin{array}{l} 
(f_1 + f_2)(x) \stackrel{{\rm def}}{=} f_1(x) + f_2(x) \; , \\
(f_1 f_2)(x) \stackrel{{\rm def}}{=} f_1(x) f_2(x)  \; , \\
(\alpha f_2)(x) \stackrel{{\rm def}}{=} \alpha(x)  f_1(x) \; . \end{array}\] 
Here $f_1, f_2 \in RMaps_{G}(X, A)$, $\alpha \in K$, $x$ is an element of $X$ 
in general position, and the operations in the right hand sides of the above
formulas are performed in $A$.  The correspondence $X \mapsto
RMaps_G(X, V)$ completes the triangle in the above
diagram; see~\cite[Proposition 8.6 and Lemma 12.3]{re}.

\begin{example} \label{ex.pgl_n}
$A = \Mn(k)$ is the algebra of $n \times n$-matrices over $k$. 
By a theorem of Wedderburn, $B$ is a $K$-algebra of type $\Mn$ 
if and only if $B$ is a central simple $K$-algebra; 
see e.g.,~\cite[Theorem 1.1]{boi} or~\cite[X.5, Proposition 7]{serrelf}.

Note that $\Aut_k(A) = \PGLn$. Thus, if $K$ is a finitely generated 
field extension of $k$, every central simple $K$-algebra $B$ of degree $n$
is of the form $RMaps_{\PGLn}(X, \Mn(k))$ for some generically 
free irreducible $\PGLn$-variety $X$ with $k(X)^G = K$. 
Moreover, the $G$-variety $X$ is uniquely determined up to birational 
isomorphism. For example, if
$B = \UD(n)$ is the universal division algebra of degree $n$ 
generated by two generic matrices then $B = RMaps_{\PGLn}(X, \Mn(k))$,
where $X = \Mn(k) \times \Mn(k)$ and $\PGLn$ acts on $X$ by 
simultaneous conjugation. This description of $\UD(n)$ is due to
Procesi; see~\cite[Theorem 2.1]{procesi}. 
\end{example} 

\begin{example} \label{ex.g-2}
$A = {\Oct}$ is the 8-dimensional split Cayley algebra
(otherwise known as the split octonion algebra). Then
$\Aut(\Oct)$ is the exceptional group $G_2$.  
By a theorem of Zorn, $B$ is a $K$-algebra of type $\Oct$ 
if and only if $B$ is a Cayley algebra over $K$. 
Cayley algebras are thus in natural 1---1 correspondence with generically
free irreducible $G_2$-varieties; see~\cite[Remark 11.4]{re}, 
\cite[Section~8.1]{serrepp},~\cite[Proposition~33.24]{boi}. 
\end{example}

\begin{example} \label{ex.f-4}
$A$ is the 27-dimensional (split) Albert algebra
(otherwise known as an exceptional simple Jordan algebra) defined over $k$.  
Then $\Aut(A)$ is the exceptional group $F_4$. Algebras of type $A$ 
are precisely
the Albert algebras, i.e., a 27-dimensional exceptional simple
Jordan algebras; see e.g.,~\cite[p. 517]{boi}, \cite[Section 9]{serrepp}.
\end{example}

\begin{remark} The results of this section remain valid if the algebra $A$ 
is replaced by a more general algebraic object consisting of a vector space 
with a tensor on it. Such objects are called structured spaces in~\cite{re}.
We refer the reader there for details; see also~\cite[Section 29]{boi}.
\end{remark}

\section{Splitting fields}
\label{sect:splt-fields}

\begin{defn} \label{def1.1}
Let $G$ be an algebraic group and $X$ be a primitive generically free
$G$-variety, $K = k(X)^G$ and $cl(X)$ = the class of $X$ in 
$H^1(K, G)$.
We will say that $X$ is {\em split} if $X$ is birationally isomorphic 
to $X/G \times G$ (as a $G$-variety). Equivalently, $X$ is split
if there exists a rational section $X/G \brokrarr X$ or, 
if $cl(X) = 1$; see, e.g., \cite[1.4.1]{popov}.

A field extension $L$ of $K$ is called a {\em splitting field} of $X$
if $X_L$ is split. Equivalently, $L$ is a splitting field 
if the image of $cl(X)$ under the natural map $H^1(K, G) \lra H^1(L, G)$
is trivial. 
\end{defn} 

\begin{remark} \label{rem.special}
Recall that an algebraic group $G$ is called {\em special} if 
every generically free $G$-variety is split. Special groups
were studied by Serre~\cite{serre1} and
classified by Grothendieck~\cite[Theorem 3]{grothendieck}; see also
\cite[Theorem 2.8]{pv}. In particular, 
$\GL_n$, $\SL_n$, $\Sympl_{2n}$, and the additive group 
${\rm G_a}$ are special. Moreover, it is easy to see that if
$N$ is a normal subgroup of $G$ and both $N$, $G/N$ are special,
then so is $G$. In particular, every connected solvable group 
is special; cf.~\cite[Section 2.6]{pv}. 
\end{remark}

\subsection*{Split $G$-varieties}

\begin{lem} \label{lem1.2}
Let $G$ be an algebraic group, let $X$ be a split G-variety and
let $H$ be a Levi-commutative subgroup 
of $G$ (see Definition~\ref{def.levi}).
If $H$ has a smooth fixed point in $X$ then

\smallskip
(a) $H$ is contained in a Borel subgroup of $G$. Moreover,

\smallskip
(b) if $H$ is diagonalizable then it is contained in a maximal 
torus of $G$.
\end{lem}

\begin{pf} 
Since $X$ is split, it is birationally isomorphic to $X_0 \times G$,
where $X_0$ is a projective variety (a projective model for $X/G$).
Let $B$ is a Borel subgroup of $G$. Consider the rational 
$G$-equivariant map
\[ X \simeq X_0 \times G \brokrarr X_0 \times G/B \; . \]
Since $X_0 \times G/B$ is a complete variety,
the Going Down Theorem~\ref{going-down} tells us that
$H$ fixes a point of $X_0 \times G/B$. Consequently,
$H$ is contained is a conjugate of $B$, and part (a) follows. 
Part (b) is immediate from part (a) 
and~\cite[Theorem~10.6(5)]{borelbook}. 
\end{pf}

\begin{remark} \label{rem.split} Note that our proof of
Lemma~\ref{lem1.2} relies only on the Going 
Down Theorem and not on Proposition~\ref{prop.proj}. 
In particular, 
Lemma~\ref{lem1.2} is valid in arbitrary characteristic;
see Remark~\ref{rem.char}.
\end{remark}

\subsection*{Proof of Theorem~\ref{thm1.2a}}
Consider the natural projection map $X_L \brokrarr X$ of $G$-varieties.
Here $L$ is a splitting field for $X$, i.e., $X_L$ is split.
By Proposition~\ref{prop.proj},
we may assume without loss of generality that $X_L$ 
is smooth and projective. Let $p^e$ be the maximal 
power of $p$ dividing $[L:K] = \deg(X_L/X)$. 
By the Going Up Theorem~\ref{going-up}(a), there exists 
a point $y \in X_L$ and a subgroup $H'$ of $H$ such that
(i) $[H:H']$ divides $p^e$ and (ii) $H' \subset \Stab(y)$. Now (i) says that
$[H:H']$ divides $[L:K]$ and (ii), in combination with Lemma~\ref{lem1.2}(b),
says that $H'$ is toral. 
\qed

\subsection*{Nontoral subgroups} 

\begin{defn} \label{def.depth} Let $G$ be an algebraic group and let 
$H$ be an abelian $p$-subgroup of $G$. The {\em depth\/} of $H$ is defined as
the smallest integer $i$ such that $H$ has a toral subgroup of index $p^i$.

Recall that a prime number $p$ is called a {\em torsion prime} for $G$
if $G$ has a nontoral abelian $p$-subgroup $H$, i.e., an abelian $p$-subgroup
$H$ of depth $\geq 1$. Following \cite[1.3]{serresubgr} we will
denote the set of torsion primes for $G$ by $Tors(G)$.
\end{defn}

\begin{remark} \label{rem.torsion}
Torsion primes have been extensively studied; see~\cite{borel},
\cite{ss}, \cite{steinberg}, \cite{serresubgr}. 
In particular, $Tors(G) = Tors(G')$ if $G'$ is a
derived subgroup of $G$; see~\cite[1.3.2]{serresubgr}. If $f \colon
\overline{G} \lra G$ is the universal cover of $G$ then $Tors(G)$ is the union
of $Tors(\overline{G})$ and the set of prime divisors of $\Ker(f)$;
see~\cite[1.3.3]{serresubgr}. 

For simply connected simple groups 
the torsion primes are given by the following table:
\[  \begin{array}{lcc}
\underline{\text{Type of $G$}} & & \underline{Tors(G)}  \\
      & &       \\
\text{$A_n$ or $C_n$} &  & \emptyset  \\
      & &       \\
B_n \; (n \geq 3) \, , \; D_n \; (n \geq 4) \, \; \text{or} \; G_2 &  & 
\{ 2 \}  \\
      & &        \\
F_4, \, E_6, \; \text{or} \; E_7 & & \{ 2, 3 \}  \\
      & &       \\
E_8 & & \{ 2, 3, 5 \}  \end{array} \] 
For details see \cite[Proposition 4.4]{borel},~\cite[Corollary 1.13]{steinberg}
or~\cite[1.3.3]{serresubgr}. 
\end{remark}

Using the terminology of Definition~\ref{def.depth}, we can rephrase
Theorem~\ref{thm1.2a} as follows.

\begin{thm} \label{thm.torsion} 
Let $G$ be an algebraic group, $H$ be a finite abelian $p$-subgroup of  
$G$ of depth $d$, $X$ be a generically free $G$-variety and $K = k(X)^G$.
If $X$ has a smooth $H$-fixed point and $L/K$ is a splitting field 
of $X$ then $p^d \mathrel| [L:K]$.  

In particular, if $X = V$ is a generically free linear representation of
$G$ then $[L:K]$ is divisible by every torsion prime of $G$.
\qed
\end{thm}

Let $G$ be an algebraic group and let $S$ be a special subgroup containing
$G$. (For example, $S$ can be taken to be $\GL_n$, $\SL_n$ or $\Sympl_{2n}$.)
We shall view $S$ as a $G$-variety with respect to the left multiplication
action; it is easy to see that this variety is generically free.
 
\begin{cor} \label{cor.torsion}
Let $G$ be an algebraic group, $H$ be a finite abelian $p$-subgroup of  
$G$ of depth $d$ and $S$ be a special group containing $G$, as above.
Suppose $K = k(S)^G$ and $L/K$ is a splitting field for $S$ (as a $G$-variety).
Then $p^d \mathrel| [L:K]$.
\end{cor}

\begin{pf} Let $\overline{S}$ be a smooth projective model of $S$
(as an $S$-variety); see Proposition~\ref{prop.proj}.
Let $V$ be a generically free linear representation of $S$.
Since $S$ is special, $V$ is split as an $S$-variety. Thus there exists
an $S$-equivariant dominant rational map 
$f \colon V \brokrarr S \simeq \overline{S}$.
Since $G \subset S$, we can view $f \colon V \brokrarr \overline{S}$ as 
a dominant rational map of $G$-varieties. Since $V$ has a $H$-fixed point,
the Going Down Theorem~\ref{going-down} tells us that so does $\overline{S}$.
Applying Theorem~\ref{thm.torsion} to the smooth variety $X = \overline{S}$,
we conclude that $p^d \mathrel| [L:K]$, as claimed.
\end{pf}

\begin{remark} \label{rem.upper}
Let $G$ be a simple group. A theorem of Tits
asserts that every $G$-variety $X$ can be split by an extension 
$L/K$, where $K = k(X)^G$ and every prime factor 
of $[L:K]$ lies in $Tors(G)$; see~\cite[2.3]{serrepp}. 
This gives a partial converse to Theorem~\ref{thm.torsion}.

More precisely, the results of~\cite{tits} show that every $G$-variety 
can be split by an extension of degree dividing $n(G)$, where $n(G)$
is given by the following table.
\[ \begin{array}{ccccc}
\underline{\text{\rm Type}} & \underline{\text{\rm Simply Connected}} &
\underline{\text{\rm Not Simply Connected}} \\
    & & \\
A_n & 1 & n+1 \\
    & & \\
B_n & 2^{\sup(1,n-4)} & 2^n \\
    & & \\
C_n & 1    & 2^{v_2(n) + 1} \\       
    & & \\
D_n & 2^{\sup(1,n-5)} & 2^{v_2(n)+n} \\
    & & \\
G_2 & 2 & - \\
    & & \\
F_4 & 6 & - \\
    & & \\
E_6 & 6  & 2\cdot 3^4     \\
    & & \\
E_7 & 12 & 2^5 \cdot 3      \\
    & & \\
E_8 & 2^7\cdot 3^3 \cdot 5 & - 
\end{array} \]
Here $G$ is an almost simple group of the indicated type 
(recall that $G$ is almost simple if the center $Z(G)$ is finite 
and $G/Z(G)$ is simple) and
$v_2(m)$ denotes the highest power of $2$ dividing $m$.

Note that the terminology of~\cite{tits} is somewhat different from 
ours. A primitive $G$-variety $X$ corresponds to a group of inner type 
over $k(X)^G$
(and if $G$ is simply connected, then of strongly inner type). With 
these conventions, the entries for all group types other 
than $E_8$ come directly from
\cite[Proposition~A1]{tits} or from \cite[Propositions~1 and~2]{tits}.

Our entry for $E_8$ follows from~\cite[Corollaire 2]{tits} and the fact
that $n(G)$ may be taken to be the degree of a splitting field for one
particular ``generic" $E_8$-variety; see~\cite[Proposition 8]{titsann}.
(Recall that we are working in characteristic zero.)
In fact~\cite[Corollaire 2]{tits} implies that $n(E_8)$ may be taken to 
be one of the numbers $2^7\cdot3\cdot5$, $2^6\cdot3^2\cdot5$ or 
$2^4\cdot3^3\cdot5$ (it is not currently known which one). The entry
for $n(E_8)$ in our table is the least common multiple of these three 
numbers.
\end{remark}

Combining the above-mentioned results of~\cite{tits} with 
Theorem~\ref{thm.torsion}, we obtain the following upper bound on the depth
of abelian $p$-subgroups of quasi-simple algebraic groups.

\begin{cor} \label{cor:upper.depth}
Let $G$ be an almost simple algebraic group and let $H$ be a finite
abelian $p$-subgroup of $G$ of depth $d$ (not necessarily elementary). 
Then $p^d$ divides the number $n(G)$ given in Remark~\ref{rem.upper}.
\qed
\end{cor}

\section{Examples}
\label{sect.examples}

In this section we illustrate Theorem~\ref{thm1.2a} for several classes
of groups.
The application of Theorem~\ref{thm1.2a} to the case $G=\PGLn$ will
come up later, after we discuss the nontoral subgroups of $\PGLn$ in
Section~\ref{sect9}; see Lemma~\ref{lem.index} below.

\subsection*{Orthogonal groups:
splitting fields of quadratic forms}

Let $K$ be a finitely generated field extension of $k$.
Recall that quadratic forms $q$ over $K$ are in 1---1 correspondence with 
primitive generically free $O_n$-varieties $X$ such that $k(X)^{O_n} = K$;
see~\cite[III. Appendix 2.2]{serregc} or~\cite[Section 29E]{boi}. 
In particular, a field extension $L/K$ splits the form if and only 
if it splits the corresponding variety.

\begin{prop} \label{prop.o_n} 
{\em Let  $q = \lf< a_1, \dots, a_n \r>$ be a quadratic form over $K$. Then 

\smallskip
(a) there exists a splitting field $L$ for $q$ such that $L/K$ is
a Galois extension with $\Galois(L/K) = (\bbZ/2\bbZ)^l$ for some
$l \leq [\frac{n+1}{2}]$. 

\smallskip
(b) Suppose $a_1, \dots, a_n$ are algebraically independent variables
over $k$, $K = k(a_1, \dots, a_n)$ and $L/K$ is a splitting field for $q$. 
Then $[L:K]$ is divisible by $2^{[\frac{n+1}{2}]}$.} 
\end{prop}

\begin{pf} (a) Suppose $n = 2m$ is even. Let
\[ L = K(\sqrt{-a_2/a_1}, \dots, \sqrt{-a_{2m}/a_{2m-1}}) \, .   \] 
Note that $L$ is a Galois extension of $K$ and
$\Galois(L/K) = (\bbZ/2\bbZ)^l$ for some $l \leq m$.  Moreover, 
since $\lf< a_{2i-1}, a_{2i} \r> \simeq \lf< a_{2i-1}, -a_{2i-1} \r>$
for $i=1,\dots,m$,
\[ q \simeq 
\lf< a_1, - a_1 \r> \oplus \dots \oplus 
\lf< a_{2m-1}, - a_{2m-1} \r> \simeq 0 \]
in the Witt group of $L$. This shows that $q$ splits over $L$.

If $n = 2m + 1$ is odd then a similar argument shows that
\[ L = k(\sqrt{-a_2/a_1}, \dots, \sqrt{-a_{2m}/a_{2m-1}}, 
\sqrt{a_{2m+1}\mathstrut}) \] is a splitting field for $q$ 
with $\Galois(L/K) = (\bbZ/2\bbZ)^{l}$ and $l \leq m+1$.

\smallskip
(b) Consider the $O_n$-variety $X = \Mn$, with $O_n$ acting (linearly) on $X$ 
by left multiplication. Recall that $k(X)^{O_n} = k(b_{ij} \mid
1 \leq i \leq j \leq n)$, where $b_{ij}(x)$ 
is the dot product of the $i$th and the $j$th columns 
of the matrix $x \in \Mn$; 
see, e.g.,~\cite[Section~2.10]{dc} or~\cite[Lemma~6.4]{re}. Note that 
since $\dim(X/O_n) = \dim(X) - \dim(O_n) = \frac{n(n+1)}{2}$, 
the generators $b_{ij}$ are algebraically indepenedent over $k$. 
The quadratic form corresponding to $X$ is the ``generic form" 
$q_X = \sum_{i \leq j} b_{ij}x_i x_j$ defined over $k(b_{ij})$.
Let $Y$ be the subvariety of $X$
consisting of $n \times n$-matrices with mutually orthogonal columns. 
Then $Y$ is irreducible (see \cite[Example 3.10]{re}) and 
the quadratic form corresponding to $Y$ is the ``generic diagonal" form
$q = \sum_{i = 1}^n a_i x_i^2 $ which appears in the statement of part (b). 
Here $a_i = b_{ii}$ and $q$ is defined over $K = k(a_1,\dots,a_n)$.
We can now view the usual orthogonalization process 
in $k^n$ as an $O_n$-equivariant 
rational map $f \colon X \brokrarr Y$. That is, we view a matrix 
$x \in \Mn$ as a collection of $n$ column vectors. 
To construct $f(x) \in Y$, we apply the orthogonalization process 
to this collection; the resulting $n$ mutually orthogonal vectors 
form the columns of $f(x)$ (see~\cite[Example 3.10]{re} for details). 

Note that the point $0_{n \times n}$ is a smooth point of $X$ fixed by
all of $O_n$; here $0_{n \times n}$ is the $n \times n$ zero matrix. 
Let $Y'$ be a smooth projective model of $Y$ 
(see Proposition~\ref{prop.proj}); we can thus think of $f$ 
as an $O_n$-equivariant rational map $X \brokrarr Y'$.
Let $H \simeq (\bbZ/2\bbZ)^n$ be the diagonal subgroup of $O_n$. 
By the Going Down Theorem~\ref{going-down}, $Y'$ has an $H$-fixed point. 
This point is smooth because every point of $Y'$ is smooth. Since $L$ 
splits $q$, it splits the $O_n$-variety $Y$ or, equivalently, the 
$O_n$-variety $Y'$. Thus Theorem~\ref{thm.torsion} tells 
us that $2^d \mathrel| [L:K]$, where $d$ is the depth of $H$.
Since the dimension of any maximal torus is $O_n$ is $[\frac{n}{2}]$, 
$d \geq n - [\frac{n}{2}] = [\frac{n+1}{2}]$,
as claimed.  
\end{pf}

The same argument with the group $SO_n$ in place of $O_n$ yields 
the following variant of Proposition~\ref{prop.o_n}. Note that 
elements $H^1(K, SO_n)$ represent equivalence classes of quadratic forms
of determinant 1; cf.~\cite[Example 8.4(b)]{re} or~\cite[(29.29)]{boi}. 

\begin{prop} \label{prop.so_n} 
{\em Let $K = k(a_1, \dots, a_n)$, $q = \lf< a_1, \dots, a_n \r>$ 
be a quadratic form of determinant 1 over $K$. Then 

\smallskip
(a) there exists a splitting field $L$ for $q$ such that $L/K$ is
a Galois extension with $\Galois(L/K) = (\bbZ/2\bbZ)^l$ for some
$l \leq [\frac{n-1}{2}]$. 

\smallskip
(b) Suppose $a_1, \dots, a_{n-1}$ are algebraically independent variables
over $k$, $a_n = (a_1 \dots a_{n-1})^{-1}$, and $L/K$ is a splitting 
field for $q$. Then 
$[L:K] \geq 2^{[\frac{n-1}{2}]}$.} 
\qed
\end{prop}

\subsection*{Exceptional groups $G_2$, $F_4$, $3E_6$ and $2E_7$}

Let $V$ be a generically free linear representation of $G$, let
$K = k(V)^G$ and let $L/K$ be a splitting field of $V$.
Theorem~\ref{thm.torsion} tells us that $[L:K]$ is divisible 
by $2$ if $G= G_2$ and by $6$ if $G = F_4$, $3E_6$ or $2E_7$.  
(Here $3E_6$ and $2E_7$ denote the simply connected groups of 
type $E_6$ and $E_7$ respectively.) If $G = G_2$, $F_4$ or $3E_6$
then this result is sharp. In other words, $V$ can be split by $L/K$
such that $[L:K]$ equals 2, if $G = G_2$ and 6, if $G = F_4$ or $3E_6$;
see Remark~\ref{rem.upper}.

\subsection*{Exceptional group $E_8$}

Recall that by a theorem of Adams~\cite{adams} $E_8$ has two maximal
elementary abelian 2-subgroups (up to conjugacy): $D(T)$ of rank 9 and $EC^8$
of rank 8. Here we are following the notational conventions 
of~\cite[Section 2]{adams}; in particular, $D(T)$ means 
``double 2-torus" and $EC^8$ means ``exotic candidate of rank 8".  
By construction $D(T)$ has depth 1. 

\begin{prop} \label{prop.ec^8} The subgroup $EC^8 \simeq (\bbZ/2\bbZ)^8$ 
of $E_8$ has depth $2$.
\end{prop}
 
Our proof of this proposition uses the theory of quadratic forms 
over $\bbZ/2\bbZ$.  Recall that if $q$ is 
a quadratic form on $V = (\bbZ/2\bbZ)^m$, the associated symmetric
(or, equivalently, skew-symmetric)
bilinear form $b_q \colon V \times V \lra V$ is defined by $b_q(v, w)
= q(v + w) - q(v) - q(w)$. Note that the ``usual" relationship between
$q$ and $b_q$ breaks down in characteristic 2: 
in particular, $b_q(v, v) = 0$ for any $v \in V$. 
The kernel of $b_q$ is called {\em the radical} of $q$.

We shall say that $v \in V$ is an {\em anisotropic} vector for $q$ if
$q(v) = 1$ and an {\em isotropic} vector if $q(v) = 0$. 
In the sequel we shall be interested in counting 
the number of anisotropic vectors for a given form $q$.
This is not a difficult task (at least in principle) because 
$q$ can always be written as a direct sum of
quadratic forms of dimension 1 and 2
(see, e.g.,~\cite[Theorem 1.4.3]{pfister}), and
if $q = r \oplus s$ then a simple counting argument shows that 
\begin{equation} \label{e.count}
\begin{array}{l}
|q^{-1}(0)| = |r^{-1}(0)| \cdot |s^{-1}(0)| + |r^{-1}(1)| \cdot |s^{-1}(1)| \\ 
    \\
|q^{-1}(1)| = |r^{-1}(1)| \cdot |s^{-1}(0)| + |r^{-1}(0)| \cdot |s^{-1}(1)|  
\; . \end{array}
\end{equation}

\begin{lem} \label{lem.quad}
Let $q$ be a quadratic form on $V = (\bbZ/2\bbZ)^7$. Suppose 
the radical of $q$ has dimension 1. 
Then $q$ has $56$, $64$ or $72$ anisotropic vectors in $V$.
\end{lem}

\begin{pf} Write $q \simeq q_1 \oplus q_2 \oplus q_3 
\oplus \lf<e \r>$, where $q_1$, $q_2$ and $q_3$ are regular 2-dimensional
quadratic forms and $e = 0$ or $1$; see, e.g.,~\cite[Theorem 1.4.3]{pfister}.
 (Here, the 1-dimensional form
$\lf<e\r>$ is the radical of $q$.) Note that over $\bbZ/2\bbZ$ there are 
only two classes of regular 2-dimensional quadratic forms: the hyperbolic
form $h$ given by $h(x_1, x_2) = x_1 x_2$ and the anisotropic form
$a(x_1, x_2) = x_1^2 + x_2^2 + x_1 x_2$. Since $a \oplus a \simeq h \oplus h$
(see~\cite[Example 2.4.5]{pfister}), 
we may assume without loss of generality that $q_2 = q_3 = h$.

\smallskip
Case 1: $e = 1$. Using \eqref{e.count} it is easy to see that
if $q_0$ is any quadratic form on
$(\bbZ/2\bbZ)^m$ then $q = q_0 \oplus \lf< 1 \r>$ has exactly $2^m$ 
anisotropic vectors in $(\bbZ/2\bbZ)^{m+1}$. 
In our situation $q_0 = q_1 \oplus q_2 \oplus q_3$ 
and $m = 6$; thus we conclude that $|q^{-1}(1)| = 64$. From now on we shall 
assume that $e = 0$. 

\smallskip
Case 2. $q = h \oplus h \oplus h \oplus \lf< 0 \r>$. 
We apply \eqref{e.count} to this form
recursively. Since  $|h^{-1}(0)| = 3$ and $|h^{-1}(1)| = 1$, we obtain 
$|q^{-1}(1)| = 56$.

\smallskip
Case 3. $q = a \oplus h \oplus h \oplus \lf< 0 \r>$. 
We note that $|a^{-1}(0)| = 1$ 
and $|a^{-1}(1)| = 3$, and apply~\eqref{e.count} recursively, to conclude 
that $|q^{-1}(1)| = 72$. 

\smallskip
This completes the proof of the lemma.
\end{pf}

\begin{pf*}{Proof of Proposition~\ref{prop.ec^8}}
Recall that $EC^8 = A_1 \times A_2 \subset G_2 \times F_4 \subset E_8$, where
$A_1 = (\bbZ/2\bbZ)^3$ is the unique (up to conjugacy) nontoral abelian
$2$-subgroup of $G_2$ and $A_2$ is the unique 
(again, up to conjugacy) nontoral abelian $2$-subgroup of $F_4$;
see~\cite[Theorem 2.17]{griess}. 
Thus, $A_1$ has a subgroup of index 2 which is toral in
$G_2$, and $A_2$ has a subgroup of index 2 which is toral in $F_4$. Taking
a direct product of these toral subgroups, we construct a subgroup of $EC^8$
of index 4 which is toral in $G_ 2 \times F_4$ and, hence, in $E_8$.
This proves that the depth of $EC^8$ is $\leq 2$.

It remains to show that the depth of $EC^8$ is $\geq 2$. Recall that
elements of $E_8$ of order 2 fall into two conjugacy classes: class $A$ 
and class $B$; cf.~\cite[Section 5]{adams} or~\cite[(2.14)]{griess}. 
If $T$ is a maximal torus in $E_8$ and $T_{(2)} =
\{ t \in T \, : \, t^2 = 1 \}$ then we have a naturally defined
$W_{E_8}$-invariant quadratic form $q$ on $T_{(2)} \simeq (\bbZ/2\bbZ)^8$; 
see~\cite[Definition 2.15]{griess}. By~\cite[Lemma 2.16]{griess} 
this form is non-singular and has maximal Witt index; moreover, 
an element $x$ of $T_{(2)}$ is of type $A$ in $E_8$ if $q(x) = 1$
and of type $B$ if $q(x) = 0$. 

In particular, of the 255 nonidentity elements of $T_{(2)}$, 120 are of
type $A$ and 135 are of type $B$. On the other hand, of the 
255 nonidentity elements of $EC^8$, 56 are of type $A$ and 199 are 
of type $B$; see~\cite[Section 5]{adams}.

We now proceed to prove that the depth of $EC^8$ is $\geq 2$. Assume, to
the contrary, that $EC^8$ has a toral subgroup $U$ of rank 7. 
Since $q$ is non-singular (i.e., the associated symplectic form $b_q$ is
non-degenerate) on $T_{(2)}$, the radical of $q_{|U}$ is of 
dimension $\leq 1$. On the other hand, since $\dim(U)$ is odd, the radical
of $q_{|U}$ cannot be trivial; thus it has dimension exactly 1. 
By Lemma~\ref{lem.quad}, $q$ has at least 56 anisotropic vectors in $U$, i.e., 
$U$ has at least 56 elements of type $A$. On the other hand, $EC^8$ 
has exactly 56 elements of type $A$. We therefore conclude that
every element of type $A$ lies in $U$. We claim that this is 
impossible because the elements of type $A$ generate $EC^8$. 
This contradiction will complete the proof of the proposition.

To prove the claim, recall that $EC^8 = A_1 \times A_2$, where
$A_1 \simeq (\bbZ/2\bbZ)^3$ lies in $G_2$ and $A_2 = (\bbZ/2\bbZ)^5$
lies in $F_4$, as above. Moreover, $A_2$ has a subgroup $R$ of order 4 
(called the radical of $EC^8$) such that 
\[ S = (A_2 - R) \cup (A_1R - R)  \]
is precisely the set of elements of $EC^8$ of type $A$; 
see~\cite[Theorem 2.17]{griess}.  We want to show 
that $\lf< S \r> = EC^8$. Indeed, $A_2 - R$ contains 28 of the 32
elements of $A_2$; these elements clearly generate all of $A_2$.
Thus $A_2 \subset \lf<S \r>$. In particular, $R \subset \lf<S \r>$. 
Now $R$, together with $A_1R$ generate $A_1$. We thus conclude that 
both $A_1$ and $A_2$ lie in $\lf<S \r>$. This proves 
that $EC^8 = A_1 \times A_2 = \lf<S \r>$, as claimed.
\end{pf*}

We are now ready to give an alternative proof of a theorem of Serre.

\begin{cor} (Serre, see \cite[Proposition~9, p.~30]{titsann} 
or~\cite[p.~1132]{tits})
\label{cor.e-8}
Suppose $E_8 \hookrightarrow S$, where $S = \GL_n$, $\SL_n$ or $\Sympl_{2n}$
for some $n$.  We shall view $S$ as an $E_8$-variety via the
left multiplication action. Suppose $K = k(S)^{E_8}$ and $L/K$ is a splitting
field of $S$. Then $[L:K]$ is divisible by 60.
\end{cor}

\begin{pf} Recall $2$ and $3$ are torsion primes of $E_8$, i.e., 
$E_8$ has an abelian $3$-subgroup and an abelian
$5$-subgroup, both of depth $\geq 1$.
Moreover, by Proposition~\ref{prop.ec^8}, $E_8$ contains an abelian
2-subgroup of depth 2.
Thus Corollary~\ref{cor.torsion} tells us that $[L:K]$ is divisible by 
$2^2 \cdot 3 \cdot 5 = 60$.
\end{pf}

\subsection*{Exceptional group $E_7$ (adjoint)}

We will now show that the (adjoint) group $E_7$ has an elementary abelian
2-subgroup of depth $\geq 2$. We begin with the following lemma.

\begin{lem} \label{lem.factor}
Let $f \colon G \lra G'$ be a surjective homomorphism of algebraic groups, such that
$\Ker(f)$ is special. Suppose $H$ is a finite abelian subgroup of $G$ and
$f(H)$ is toral in $G'$. Then $H$ is toral in $G$. 
\end{lem}

\begin{pf} Suppose $f(H) \subset T' \subset G'$, where $T'$ be a torus in $G'$.
Denote $f^{-1}(T')$ by $S$. Then $H \subset S \subset G$. Moreover,
since both $\Ker(f)$ and $S/Ker(f) \simeq T'$ are special, 
we conclude that $S$ is
special as well; see Remark~\ref{rem.special}.  This means that $H$ 
is toral in $S$ (see e.g.,~\cite[1.5.1]{serresubgr} 
or Example~\ref{ex.special}); 
hence, $H$ is toral in $G$.  \end{pf}

\begin{prop} \label{prop.e-7} The (adjoint) group $E_7$ has an elementary 
abelian 2-subgroup of depth $\geq 2$.
\end{prop}

Our proof uses the idea of Adams (see~\cite[Introduction]{adams})
to study nontoral 2-subgroups in groups of type $E_7$ by embedding $2E_7$
into $E_8$.

\begin{pf}  Let $EC^8$ be a maximal elementary abelian subgroup of $E_8$ of
rank 8, as in Proposition~\ref{prop.ec^8}. As we mentioned in the proof of that
proposition, $EC^8$ has 56 elements of type $A$ (in $E_8$). 
Let $x$ be one of these 56 elements.  Denote the centralizer
$C_{E_8}(x)$ by $C$. Note that $EC^8 \subset C$. Moreover,  
$C \simeq 2A_1 E_7$; see~\cite[p.~280]{griess}. Thus there is an exact 
sequence
\[ \{ 1 \} \lra \SL_2 \lra C \stackrel{f}{\lra} E_7 \lra \{ 1 \} \; . \]
We claim that $f(EC^8)$ has depth $\geq 2$ in $E_7$. Indeed, assume 
the contrary. Then $f(EC^8)$ contains a subgroup $H'$ of index 2 which is
toral in $E_7$.  By Lemma~\ref{lem.factor}, $H = f^{-1}(H') \cap EC^8$ is
toral in $C$ and thus in $E_8$. Since $H$ is a toral subgroup of index 2 in
$EC^8$, this implies that $EC^8$ has depth $\leq 1$, contradicting 
Proposition~\ref{prop.ec^8}.
\end{pf}

We can now prove an analogue of Corollary~\ref{cor.e-8} for $E_7$.

\begin{cor} \label{cor.e-7}
Suppose $E_7 \hookrightarrow S$, where $S = \GL_{n}$, $\SL_n$ or $\Sympl_{2n}$
for some $N$.  We shall view $S$ as an $E_7$-variety via the
left multiplication action. Suppose $K = k(S)^{E_7}$ and $L/K$ is a splitting
field of $S$. Then $[L:K]$ is divisible by 12.
\end{cor}

\begin{pf} Recall $E_7$ has a nontral abelian $3$-subgroup, 
i.e., a $3$-subgroup of depth $\geq 1$; see, e.g.,~\cite{griess}.
Moreover, by Proposition~\ref{prop.e-7}, $E_7$ contains an abelian
2-subgroup of depth $\geq 2$.  Thus Corollary~\ref{cor.torsion} tells 
us that $[L:K]$ is divisible by $2^2 \cdot 3 = 12$.
\end{pf}

\section{Splitting groups}
\label{sect:splt-groups}

\subsection*{Definition and first examples}

\begin{defn} \label{def1.5} 
Let $X$ be a generically free primitive $G$-variety and
let $K = k(X)^G$. We shall say that a finite group $A$ is a {\em splitting
group} for $X$ if there exists a splitting field $L$ for $X$ such that
$L/K$ is (finite) Galois and $\Galois(L/K) = A$.
\end{defn}

\begin{example} \label{ex.finite} {\em Let $G$ be a finite group
and let $X$ be a generically free irreducible $G$-variety. 
Then $G$ is a splitting group for $X$.}
\end{example}
 
\begin{pf}
Suppose $L = k(X)$ and $K = k(X)^G$. Then $k(X_L) = L \otimes_K L =
L \oplus \dots \oplus L$ ($|G|=[L:K]$ times). In other words, 
up to birational equivalence, $X_L$ is the disjoint
union of $|G|$ copies of $Y$ and $G$ acts on $X_L$ by permuting these
copies.  Consequently, $X_L$ is split 
as a $G$-variety and $G = \Galois(L/K)$ is a splitting group.
\end{pf}

\begin{example} \label{ex.galitskii}
{\em Let $G$ be a (connected) semisimple group, and let $W$ be the Weyl 
group of $G$. Then every irreducible generically free 
$G$-variety $X$ has a splitting group
which is isomorphic to a subgroup of $W$.}
\end{example}

\begin{pf} \label{rem.gal} 
Let $X$ be a generically free irreducible $G$-variety
and let $ \pi \colon X \brokrarr X/G $ be the rational quotient map.
An irreducible subvariety $S$ of $X$ is called 
a {\em Galois section} if $GS$ is dense in $X$, 
i.e., $\pi_{|S}$ is dominant,
and the field extension $k(S)/k(X)^G$ induced by $\pi$, is a finite
Galois extension. We shall denote the group 
$\Galois(k(S)/k(X)^G)$ by $\Galois(S)$.  

A theorem of Galitskii asserts that every $G$-variety $X$ has a Galois section
$S$; see~\cite{galitskii}. Moreover, by \cite[Remark 1.6.3]{popov}
$S$ can be chosen so that $\Galois(S)$ is isomorphic to
a subgroup $H$ of $W$.
It is easy to see that in order to split $X$ as a $G$-variety it is
sufficient to split $S$ as a $\Galois(S)$-variety. Now
Example~\ref{ex.finite} tells us that $H$ is a splitting group for $X$.
\end{pf}

If $G$ is a connected but not necessarily semisimple then the assertion of
Example~\ref{ex.galitskii} remains true if we define $W$ as the Weyl
group of $G_{ss} = G/R(G)$, where $R(G)$ is the radical of $G$.

\subsection*{Two elementary lemmas from group theory}

Our next goal is to prove Theorem~\ref{thm1.2b}. We begin with two 
elementary lemmas.

\begin{lem} \label{lem.sub=quot}
Let $P$ be a finite abelian group.

\smallskip
(a) Every quotient group of $P$ is isomorphic to a subgroup of $P$.

\smallskip
(b) Every subgroup of $P$ is isomorphic to a quotient group of $P$.
\end{lem}

\begin{pf} Suppose $Q$ is a quotient group of $P$. Then every character
of $Q$ lifts to a character of $P$. This gives an inclusion $Q^{\ast}
\hookrightarrow P^{\ast}$ of dual groups. Since $Q^{\ast} \simeq Q$ and
$P^{\ast} \simeq P$, part (a) follows. Part (b) is proved in a similar
manner.
\end{pf}

\begin{lem} \label{lem2.1} Let $A$ and $B$ be (abstract) groups and
$A \times B$ acts on a set $Z$. Let $W$ be the set of $A$-orbits in $Z$
and let $f \colon Z \lra W$ be the natural projection.
Assume $z \in Z$ and $w = f(z)$.  Then 

\smallskip
(a) $\Stab_B(z)$ is a normal subgroup of $\Stab_B(w)$. 

\smallskip
(b) Let $S = \Stab_{A \times B}(z)$. Then $\Stab_B(z)$ is normal in $S$ and
$S/\Stab_B(z)$ is isomorphic to a subgroup of $A$; we shall denote this 
subgroup by $A_0$.

\smallskip
(c) $\Stab_B(w)/\Stab_B(z)$ is isomorphic to a quotient of $A_0$. 
\end{lem}

\begin{pf} (a) Suppose $b \in \Stab_B(w)$. 
Since the actions of $A$ and $B$ 
on $Z$ commute, $f(bz) = bf(z) = bw = w$. Consequently, $bz = az$
for some $a \in A$ and thus
\[ b \Stab_B(z) b^{-1} = \Stab_B(bz) = \Stab_B(az) = \Stab_B(z) \, .\] 
This proves part (a).

Let $\pi_A$ and $\pi_B$ be, respectively, the natural projections 
$A \times B \lra A$ and $A \times B \lra B$. 

\smallskip
(b) The kernel of the map $\pi_A \colon S \lra A$ is $S \cap B = \Stab_B(z)$,
and part (b) follows.

\smallskip
(c) Note that $b \in \Stab_B(w)$ if and only if $bz = az$ for some
$a \in A$ or, equivalently, if $(a^{-1}, b) \in S$ for some $a \in A$.
In other words,
$\Stab_B(w) = \pi_B(S)$. Consequently, we have a surjective homomorphism
\[ \pi_B \colon A_0 = S/\Stab_B(z) \lra \Stab_B(w)/\Stab_B(z) \, . \]
This completes the proof of part (c).
\end{pf}

\subsection*{Proof of Theorem~\ref{thm1.2b}}
Let $K = k(X)^G$ and let $L/K$ be a Galois extension such that
$\Galois(L/K) = A$ and $X_L$ is split.
Note that $A \times G$ acts rationally on $X_L$. By
a theorem of Rosenlicht (see~\cite[Theorem 1]{rosenlicht1}), we can 
choose a birational model for $X_L$ so that this action becomes regular.
Moreover, after applying Proposition~\ref{prop.proj} and 
Theorem~\ref{thm.resol} to $X_L$, we may assume that

\smallskip
(i) $X_L$ is smooth and projective, and

\smallskip
(ii) for every $z \in X_L$, $\Stab_{A \times G}(z)$ is Levi-commutative
(see Definition~\ref{def.levi}).

\smallskip
Note that by our construction the map $h \colon X_L \brokrarr X$ 
is a rational quotient map for the $A$-action on $X_L$. 
Since $A$ is a finite group and $Z$ is projective, 
there exists a geometric quotient map $f \colon X_L \lra W = X_L//A$ for the
$A$-action on $X_L$ with $W$ projective; see Lemma~\ref{lem.quasi-proj}.
Note that by the universal property of categorical (and, hence, geometric)
quotients, the $G$-action on $X_L$ descends to $W$; by our construction,
$W$ and $X$ are birationally isomorphic as $G$-varieties.
Applying the Going Down Theorem~\ref{going-down}
to the birational isomorphism $X \stackrel{\simeq}{\brokrarr} W$,
we conclude that $W$ has a $H$-fixed point. Denote this point by $w$. 

Now choose $z \in f^{-1}(w) \in X_L$ and apply Lemma~\ref{lem2.1} with
$Z = X_L$ where we view $Z$ as an $A \times H$-variety via 
the obvious inclusion of $A \times H$ in $A \times G$.
By Lemma~\ref{lem2.1}(b), $A$ has a subgroup
$A_0 \simeq S/\Stab_H(w)$, where $S = \Stab_{A \times H}(z)$.
Note that $S$ is a finite subgroup of $\Stab_{A \times G}(z)$, 
and $\Stab_{A \times G}(z)$ is Levi-commutative by our construction
(see condition (ii) above). We conclude that $S$ is abelian; 
see Lemma~\ref{lem.levi}.  Thus $A_0$ is also abelian.

By Lemma~\ref{lem2.1}(c), $A_0$ has a quotient of the form
$\Stab_H(w)/\Stab_H(z) = H/\Stab_H(z)$.  Denote $\Stab_H(z)$ by $H'$. 
Since $Z = X_L$ is split, Lemma~\ref{lem1.2}(b) says that $H'$ 
is toral in $G$. Thus $H/H'$ is a quotient of $A_0$, with $H'$ toral.
By Lemma~\ref{lem.sub=quot}, $A_0$ (and, hence, $A$) has a subgroup
isomorphic to $H/H'$, as claimed.
\qed

\subsection*{Examples}

\begin{example} \label{ex.special} (cf.~\cite[1.5.1]{serresubgr})
{\em Let $G$ be a special group (see Remark~\ref{rem.special}). 
Then every finite abelian subgroup of $G$ is toral.}

Indeed, let $V$ be a generically free linear representation of $G$ and
let $H$ be a finite abelian subgroup of $G$.
Since $G$ is special, $V$ is split, i.e., $A = \{ 1 \}$ is a splitting 
group for $V$. On the other hand, since the origin of $V$
is a smooth $H$-fixed point, 
there exists a toral subgroup $H_T$ of $H$ such that
$H/H_T$ is isomorphic to a subgroup of $A =  \{1 \}$. 
In other words, $H=H_T$ is toral, as claimed.
\qed
\end{example}

\begin{example} \label{ex.o_n-3} {\em Let $a_1, \dots, a_n$ be independent 
variables over $k$ and let $q = \lf< a_1, \dots, a_n \r>$ be the generic
quadratic form of dimension $n$. Then any splitting group of $q$ contains
a copy of $(\bbZ/2\bbZ)^{[\frac{n+1}{2}]}$.}

The proof is the same as in 
Proposition~\ref{prop.o_n}(b), with Theorem~\ref{thm1.2b} 
used in place of Theorem~\ref{thm1.2a}.
\qed
\end{example}

\begin{example} \label{ex.e-8} {\em Let $G = E_7$ (adjoint) or $E_8$.
Suppose $G \hookrightarrow S$, where $S = \GL_n$, $\SL_n$ or $\Sympl_{2n}$
for some $n$.  
Then any splitting group of $S$ (viewed as a $G$-variety with respect
to the left multiplication action) contains a copy of $(\bbZ/2\bbZ)^2$}.

Indeed, $G$ has an elementary abelian 2-subgroup $H$
of depth $\geq 2$; see Propositions~\ref{prop.ec^8} and~\ref{prop.e-7}.
If $X = \overline{S}$ is a smooth projective model for $S$ (as an $S$-variety)
then the argument of Corollary~\ref{cor.torsion} shows that
$H$ has a fixed point in $X$. Thus we can apply 
Theorem~\ref{thm1.2b} to $X$.
\qed
\end{example}

\section{Abelian subgroups of $\protect\PGLn$}
\label{sect9}

The rest of this paper will be devoted to applications of Theorem~\ref{thm1.2b} 
(with $G = \PGLn$) to the theory of central simple algebras. In this section 
we lay the foundation for these applications by studying finite abelian 
subgroups of $\PGLn$. 

\subsection*{Symplectic modules}

We begin by recalling the notion of a symplectic module from~\cite{ta2}.
Let $H$ be an abelian group; in the sequel we shall refer  
to such groups as $\bbZ$-modules or just modules. We will always 
assume $H$ is finite.
A skew-symmetric form on $H$ is
a skew-symmetric $\bbZ$-bilinear map $\omega \colon
H \times H \lra \bbQ/\bbZ$. If $H$ is written multiplicatively,
we will usually identify $\bbQ/\bbZ$ with the multiplicative group
of roots of unity in $k^{\ast}$.  A subgroup (or, 
equivalently, a $\bbZ$-submodule) $H'$ of $H$ is called 
{\em isotropic} if $\omega(h, h') = 0$ for every $h, h' \in H'$. 

We will say that $\omega$ is 
{\em symplectic} if it is non-degenerate, i.e., the homomorphism
$H \lra H^{\ast}$ it defines, is an isomorphism. If $\omega$ is
symplectic then a subgroup $H'\subset H$ is called {\em Lagrangian} if
it is a maximal isotropic subgroup, i.e., if $H'$ is not contained in any
other isotropic subgroup. It is easy to see that
$H'$ is Lagrangian if and only if it is isotropic and $|H'|^2 = |H|$;
cf.~\cite[Corollary~3.1]{ta2}.

\begin{lem} \label{lem.isotr}
Let $(H, \omega)$ be a symplectic module of order $n^2$.

\smallskip
(a) If $\Lambda$ is a Lagrangian submodule then $H/\Lambda \simeq \Lambda$ (as abelian groups).

\smallskip
(b) Let $H_1$ be a subgroup of $H$, and
let $I$ be an isotropic subgroup of $H_1$.  Then $H_1/I$ contains an
isomorphic copy of $I_1$, where $I_1$ is an isotropic subgroup of $H$
and $|H_1|$ divides $n|I_1|$.
\end{lem}

\begin{pf} (a) For $h \in H$, let $\chi_h \colon \Lambda \lra k^{\ast}$ be 
the character given by $\chi_h(l) = \omega(h, l)$. 
Then $h \mapsto \chi_h$
is a group homomorphism $\phi \colon H \lra \Lambda^*$. Since $\omega$ is 
non-degenerate, $\phi$ is onto. Since $\Lambda$ is Lagrangian, 
$\Ker(\phi) = \Lambda$.  Thus $H/\Lambda \simeq \Lambda^*$. 
Since $\Lambda^{\ast} \simeq \Lambda$, the lemma follows. 

(b) We may assume without loss of generality that $I$ is a maximal
isotropic subgroup of $H_1$. Indeed, let $I_{max}$ be a maximal 
isotropic subgroup of $H_1$ containing $I$.  Suppose we can find an
isotropic subgroup $I_1$ such that
$|H_1|$ divides $n|I_1|$ and $H/I_{max}$ has a subgroup isomorphic to $I_1$.
Since $H/I_{max}$ is isomorphic to a quotient, and hence
a subgroup of $H/I$ (see Lemma~\ref{lem.sub=quot}), 
the same $I_1$ will work for $I$. 

Thus we may (and will) assume that $I$ is a maximal isotropic 
subgroup of $H_1$. Let $\Lambda$ be a Lagrangian subgroup of $H$ containing
$I$. Then $\Lambda \cap H_1 = I$ and thus $H_1/I \hookrightarrow H/\Lambda 
\simeq \Lambda$; the last isomorphism is given by part (a).  Denote the
image of $H_1/I$ in $\Lambda$ by $I_1$. Then  
\[ |H_1| = |I| \cdot |H_1/I| = |I| \cdot |I_1| \]
divides $|\Lambda| \cdot |I_1| = n |I_1|$, as claimed.
\end{pf}

\begin{defn} \label{def.S1} (cf. \cite[Section 4]{ta2})
Let $A$ be an abelian group. We
define a skew-symmetric form $\omega_A$ on $A \times A^{\ast}$ by
\[ \omega_1(a_1 \oplus \chi_1, a_2 \oplus \chi_2) = \chi_1(a_2) - \chi_2(a_1)
\; . \]
\end{defn}

\begin{lem} \label{lem.sub}
(a) $(A \times A^*, \omega_A)$ is a symplectic module, 
and $A \times \{ 1 \}$ is a Lagrangian submodule.

\smallskip
(b) Moreover, every symplectic module $H$ is of the form
$(A \times A^*, \omega_A)$ for a suitable Lagrangian submodule $A$ of $H$. 

\smallskip
(c) Let $(H, \omega)$ be a symplectic module, $H \simeq (\bbZ/p\bbZ)^{2r}$. 
If $s \leq r$ then $(H, \omega)$ has a symplectic submodule 
$(H_1, \omega_{| H_1})$ of rank $2s$.
\end{lem}

\begin{pf} Parts (a) and (b) are proved in~\cite[Section 4]{ta2}. 
Proof of (c): By part (b), we can write 
$(H, \omega)$ as $(A \times A^*, \omega_A)$, where $A = (\bbZ/p\bbZ)^r$.
Let $e_1, \dots, e_r$ be an $\bbZ/p\bbZ$-basis of $A$, viewed
as an $r$-dimensional vector space
over $\bbZ/p\bbZ$. Then the module
$H_1$ spanned by $(e_i, 1)$ and $(1, e_j^{\ast})$ as $i, j$ 
range from $1$ to $s$, has the desired properties. 
\end{pf}

\subsection*{The form $\alpha_H$}

Abelian subgroups of $\PGLn$ are naturally
endowed with a skew-symmetric bilinear form. 

\begin{defn} \label{def6.5}
Let $H$ be a finite abelian subgroup of $\PGL_n$. For $a, b \in H$ define
$\alpha_H(a, b) = ABA^{-1}B^{-1}$, where $A$ and $B$ are elements 
of $\GL_n$ representing $a$ and $b$ respectively. It is easy to see
that $\alpha_H(a, b)$ does not depend of the choice of $A$ and $B$ and
$\alpha_H \colon H \times H \lra k^{\ast}$ defined this way,
is a skew-symmetric form on $H$. 
\end{defn}

\begin{lem} \label{lem.subgr}
Let $H$ be a finite abelian subgroup of $\PGL_n$. Then
the following conditions are equivalent.

\smallskip
(a) $H$ lifts to an abelian subgroup of $\SL_n$.

\smallskip
(b) $H$ is toral. 

\smallskip
(c) The skew-symmetric form $\alpha_H \colon H \times H \lra k^{\ast}$
given in Definition~\ref{def6.5} is trivial, i.e., 
$\alpha_H(a, b) = 1$ for every $a, b \in H$.
\end{lem}

\begin{pf} (a) $\Longrightarrow$ (b). Recall that every finite abelian
subgroup of $\SL_n$ can be simultaneously diagonalized and hence, is toral. 
(Alternatively, since $\SL_n$ is a special group, this follows from
Example~\ref{ex.special}.) The tori of $\PGL_n$ are precisely the
images of the tori in $\SL_n$ under the natural projection 
$\SL_n \lra \PGLn$, and part (b) follows. 

\smallskip
(b) $\Longrightarrow$ (c). Suppose $H$ is contained in a maximal torus 
$T \subset \PGL_n$ and let $S$ be the preimage of $T$ in $\SL_n$. Then
$S$ is a maximal torus of $\SL_n$. Thus any $a, b \in H$ can be lifted
to, respectively, $A, B \in S$. Since $A$ and $B$ commute, we conclude that
$\alpha_H(a, b) = ABA^{-1}B^{-1} = 1$, as claimed.

\smallskip
(c) $\Longrightarrow$ (a). Since
$\alpha_H$ is trivial, the preimage of $H$ in $\SL_n$ is a finite abelian group.
\end{pf}

\subsection*{The embedding $\phi$}
We will now show that any symplectic module $H$ can be obtained from
an abelian subgroup of $\PGLn$, as above, with
$n = \sqrt{|H|}$.  Note that by
Lemma~\ref{lem.sub}(b) we may assume 
$H = (A \times A^*, \omega_A)$ for some abelian group $A$. 

\begin{defn} \label{def.phi_A} (cf.~\cite[Definitions 8.7 and 8.10]{ry})
Let $A$ is an abelian group of order $n$.
We define the embedding
\[
\phi\colon A \times A^*\hookrightarrow\PGLn
\]
as follows.
Identify $\PGLn$ with $\PGL (V)$, where
$V = k[A] =$ the group algebra of $A$.  The group 
$A$ acts on $V$ by the regular representation $a \mapsto P_a \in \GL(V)$,
where
\[ P_a\Bigl(\sum_{b \in A} c_b b\Bigr) = \sum_{b \in A} c_b ab  \]
for any $a \in A^{\ast}$ and $c_b \in k$. 
The dual group $A^{\ast}$ acts on $V$ by the representation
$\chi \mapsto D_{\chi} \in \GL(V)$, where
\[ D_{\chi}\Bigl( \sum_{a \in A} c_a a\Bigr) = \sum_{a \in A} c_a \chi(a) a  \]
for any $\chi \in A^{\ast}$ and $c_a \in k$.
We define $\phi$ by{\normalshape
\[
\phi(a,\chi) = \text{ the image of } P_a D_\chi \text{ in }\PGL(V)\ .
\]}
\end{defn}

\begin{lem} \label{lem:reminder.phi}
Let $A$ be a finite abelian group, $a, b \in A$ and
$\chi, \mu \in A^{\ast}$. Then

\smallskip
(a) $D_{\chi} P_a = \chi(a) P_a D_{\chi}$. 

\smallskip
(b) $(P_a D_{\chi}) (P_b D_{\mu}) (P_a D_{\chi})^{-1} = \chi(b) \mu^{-1}(a) 
(P_b D_{\mu})$

\smallskip
(c) The embedding $\phi$ of Definition~\ref{def.phi_A} is a
monomorphism of groups, and $\phi(A\times A^*)$ is subgroup of $\PGLn$
isomorphic to $A\times A^*$.
\end{lem}

\begin{pf}
See~\cite[Lemmas~8.8 and~8.11(i)]{ry}.
\end{pf}

\begin{lem} \label{lem.emb}
Let $A$ be an abelian group of order $n$ and let  
$\phi \colon A \times A^* \hookrightarrow \PGLn$ be the embedding
of Definition~\ref{def.phi_A}. Then $\phi$ induces an isomorphism of
$(A \times A^*, \omega_A)$ and $(\phi(A \times A^*),\alpha)$
as modules with skew-symmetric forms,
where $\omega_A$ is as in Definition~\ref{def.S1}
and $\alpha=\alpha_{\phi(A \times A^*)}$ is as in Definition~\ref{def6.5}.
In particular, $(\phi(A \times A^*),\alpha)$ is a symplectic module.
\end{lem}

\begin{pf} Let $h_1 = (a_1, \chi_1)$ and 
$h_2 = (a_2, \chi_2) \in A \times A^*$. Then we want to show
that \[ \alpha(\phi(h_1), \phi(h_2)) = \omega_A(h_1, h_2) = 
\chi_1(a_2) \chi_2(a_1)^{-1} \; . \]
On the other hand, by definition of $\alpha$, we have
\[ \alpha(\phi(h_1), \phi(h_2)) = (P_{a_1} D_{\chi_1}) (P_{a_2} D_{\chi_2}) 
(P_{a_1} D_{\chi_1})^{-1} (P_{a_2} D_{\chi_2})^{-1} \; . \]
The desired equality now follows from Lemma~\ref{lem:reminder.phi}(b).
\end{pf}

\begin{cor} \label{cor:depth.pgln}
 Let $A$ be an abelian group of order $n = p^r$.
 Then the subgroup $H=\phi(A\times A^*)\subset\PGLn$ is of depth $r$.
\end{cor}

\begin{pf}
Let $H_T$ be any maximal (with respect to inclusion) toral subgroup of $H$.
By Lemma~\ref{lem.subgr}, $H_T$ is isotropic; as it is maximal,
it is Lagrangian.  The index $[H:H_T]=n^2/n=n = p^r$, and hence, the depth
of $H$ is $r$.
\end{pf}

If $n = p^r$ then the depth of any $p$-subgroup of $\PGLn$ is $ \leq r$.
This can be shown directly or, alternatively, derived from 
Theorem~\ref{thm1.2a}, since any central simple algebra of degree $n$
is split by a degree $n$ extension of its center.

\section{Symplectic modules and division algebras}
\label{sect10}

We now know enough about the nontoral subgroups of $\PGLn$ to
proceed with our results on division algebras.

\subsection*{When is $RMaps_{\protect\PGLn}(X,\protect\Mn)$ a
division algebra?}

We begin with an application of Theorem~\ref{thm1.2a}.

Let $X$ be a generically free irreducible $\PGLn$-variety. Recall that
$A=RMaps_{\PGL_n}(X, \Mn)$ is a central simple algebra with the center
$Z(A)=k(X)^\PGLn$; $A$ is of the form $\Mat_s(D)$, where $D$ is a
division algebra.  The degree $d$
of $D$ is called the index of $A$, and $sd=n$.  The following 
lemma relates smooth points in $X$ fixed by finite abelian subgroups
of $\PGLn$, to the index of $A$.

Let $H$ be a finite abelian subgroup of $\PGLn$.  The skew-symmetric form
$\alpha_H$ on $H$ may be singular; the quotient $H/\Ker(\alpha_H)$ is
a symplectic module, and hence, $|H/\Ker(\alpha_H)|=m^2$ for some
integer $m$.

\begin{lem} \label{lem.index}
With the notations as above, suppose that $H$ has a smooth fixed point $x\in X$.

Then the index of $A$ is divisible by $m$.
In particular, $m\mathrel|n$, and if $H=\phi_P(P\times P^*)$ where $P$
is an abelian group of order $n$ (so that $m=n$), then $A$ is a
division algebra.
\end{lem}

\begin{pf} Let $F$ be the center of $D$ (and of $A = \Mat_s(D)$), and let
$K$ be a maximal subfield of $D$. Recall that 
$[K:F] = d = \deg(D)$ = index($A$) and that $K$ is a splitting field of $A$.
By Theorem~\ref{thm1.2a}, $d$ is divisible by $|H/H_T|$ where $H_T$ is some
toral subgroup of $H$.
By Lemma~\ref{lem.subgr}, $H_T$ is isotropic in $H$; we may assume
that $H_T$ is a maximal isotropic subgroup of $H$.
Then $H_T\supset\Ker(\alpha_H)$, and the
image of $H_T$ in $H/\Ker(\alpha_H)$ is Lagrangian; it follows that
$|H/H_T|=m$, i.e., $d$ is divisible by $m$.

The equality $sd=n$ implies then that $m\mathrel|n$, and if $m=n$ then
$s=1$, i.e., $A$ is a division algebra.
\end{pf}

\subsection*{Proof of Theorem~\ref{thm1.3}} 
The following proposition is an application of Theorem~\ref{thm1.2b}.

\begin{prop} \label{prop6.8}
Let $H$ be a finite abelian subgroup of $\PGL_n$ of order $n^2=p^{2r}$,
such that $(H, \alpha_H)$ is a symplectic module (i.e., $\alpha_H$ 
is non-degenerate on $H$; see Definition~\ref{def6.5}). Suppose
$X' \brokrarr X$ is a rational cover of irreducible generically 
free $\PGLn$-varieties, $p^e$ is the largest power of $p$ dividing 
$\deg(X'/X)$, and $X$ has a smooth point fixed by $H$.  
Then any splitting 
group $A'$ for $X'$ contains an isomorphic copy of some isotropic
subgroup $I_1 \subset H$, where $|I_1| \geq p^{r-e}$. 

In particular, if $e = 0$, $A'$ contains an isomorphic copy of a
Lagrangian subgroup of $H$.
\end{prop}

\begin{pf} By the Going Up Theorem~\ref{going-up}(a), $X'$ has an
$H_1$-fixed point for some subgroup $H_1$ of $H$ of order $p^{2r-e}$.
By Theorem~\ref{thm1.2b}, $A'$ contains a copy of $H_1/I$, 
where $I$ is a toral subgroup of $H_1$.
Lemma~\ref{lem.subgr} says that $\alpha_H$ restricted to $I$ is trivial, 
i.e., $I$ is an isotropic subgroup. Thus by Lemma~\ref{lem.isotr}(b), 
$H_1/I$ contains a copy of $I_1$, where $I_1$ is an isotropic subgroup
of $H_1$ and $|H_1|$ divides $p^r \cdot|I_1|$. Since $|H_1| \geq p^{2r-e}$,
this translates into $|I_1| \geq p^{r-e}$, as claimed.
\end{pf}

We now continue with the proof of Theorem~\ref{thm1.3}.
Recall that $\UD(n) = RMaps_{\PGL_n}(X, \Mn)$, where $X = \Mn \times
\Mn$, with $\PGLn$ acting by simultaneous conjugation; 
see Example~\ref{ex.pgl_n}.  Let 
$D = \UD(n) \otimes_{Z(n)} K$, as in the statement of Theorem~\ref{thm1.3}. 
Then $D = RMaps_{\PGL_n}(X_K, \Mn)$.
Recall that we are assuming $n = p^r$, and
$p^e$ is the highest power of $p$ dividing $[K: Z(n)] = \deg(X_K/X)$.
Also recall that $A$ is a splitting group for $D$ if and only if $A$ 
is a splitting group for $X_K$ (as a $\PGLn$-variety); see 
Definition~\ref{def1.5}.

Note that $X = \Mn \times \Mn$ 
has a smooth point (namely, the origin) fixed by all of $G$.
Let $P$ be an abelian $p$-group of order $n = p^r$; then
$H = \phi_{P}(P \times P^{\ast})$ is an abelian $p$-subgroup of
$\PGL_n$. By Lemma~\ref{lem.emb}, $(H, \alpha_H)$ is a symplectic 
module. Applying Proposition~\ref{prop6.8} to $H$ and remembering that
every symplectic module of order $p^{2r}$ is isomorphic to one of the form
$\phi_P(P \times P^*)$ for some $P$ (see Lemmas~\ref{lem.sub}(b)
and~\ref{lem.emb}), we obtain the following
generalization of \cite[Corollary 7.2]{ta1} (in characteristic 0):

\begin{prop} \label{prop.lagr}
Let $Z(p^r)$ be the center of the generic division algebra $\UD(p^r)$,
let $K$ be a field extension of $Z(p^r)$ and let
$D = \UD(p^r) \otimes_{Z(p^r)} \, K$. Suppose $p^e$ is the highest
power of $p$ dividing $[K:Z(p^r)]$, where $e$ is a non-negative integer and
$e \leq r-1$. 

If $A$ is a splitting group of $D$ then for every symplectic module 
$H$ of order $p^{2r}$, there exists an isotropic submodule $I_1$ of order
$p^{r-e}$ such that $A$ contains an isomorphic copy of $I_1$.
\qed
\end{prop}

In order to finish the proof of Theorem~\ref{thm1.3} we use a comparison 
argument, as in the proof of \cite[Proof of Theorem 7.3]{ta1}.
Let
\begin{equation} \label{e.h1-h2}
H_1 = \phi_{P_1}(P_1 \times P_1^*)\text{ and }
H_2 = \phi_{P_2}(P_2 \times P_2^*)\ ,
\end{equation}
where 
\begin{equation} \label{e.p1-p2} P_1 = (\bbZ/p\bbZ)^r \; \, \text{and}
\; \, P_2 = \bbZ/p^r\bbZ \; .
\end{equation}
By Proposition~\ref{prop.lagr}, $A$ contains an isomorphic copy $I_1$
of an isotropic subgroup
of $H_1$, and an isomorphic copy $I_2$ of an isotropic subgroup of
$H_2$, such that $|I_1| = |I_2| = p^{r-e}$.
Since $H_1 \simeq (\bbZ/p\bbZ)^{2r}$ and $H_2 \simeq (\bbZ/p^r\bbZ)^2$,
$I_1 \simeq (\bbZ/p\bbZ)^{r-e}$ and $I_2$ has rank $\leq 2$.
We may assume without loss of generality that
both $I_1$ and $I_2$
are contained in the same Sylow $p$-subgroup $A_p$ of $A$. 
Since the intersection of $I_1$ and $I_2$ 
has exponent $p$ and rank $\leq 2$, we see that 
\[
|A_p| \geq |I_1 I_2| = \frac{|I_1| \, |I_2|}{|I_1 \cap I_2|} \geq  
\frac{p^{2r-2e}}{p^2} = p^{2r-2e-2} \; ,
\]
where $I_1 I_2=\{\gamma_1\gamma_2\mid\gamma_1\in I_1,\,\gamma_2\in I_2\}$.

This shows that $|A_p|$ is divisible by $p^{2r-2e-2}$ and, hence, so is
$|A|$, as claimed.
\qed

\begin{remark} \label{rem1.3} 
The only property of $X = \Mn \times \Mn$ used in the above proof is
that each of the finite abelian subgroups
$H_1 \simeq (\bbZ/p^r\bbZ)^2$ and
$H_2 \simeq (\bbZ/p\bbZ)^r$ of $\PGL_n$
has a smooth fixed point in $X$.
Thus our argument shows that Theorem~\ref{thm1.3}
remains valid if the universal division algebra $\UD(n)$ is replaced
by the algebra $U' = RMaps_{\PGL_n}(X, \Mn)$, 
where $X$ is an irreducible generically free
$\PGLn$-variety $X$ such that $H_i$ has a smooth fixed point in $X$ 
for $i = 1, 2$. There are many choices for such $X$; in particular,
Proposition~\ref{prop.trdeg6} shows that $X$ can be chosen so that
$\dim(X/\PGLn) = 2r$ or, equivalently, $\trdeg_k(Z(U')) = 2r$, where
$Z(U')$ is the center of $U'$.
\end{remark}

\begin{remark} \label{rem:thm1.4-ext}
Tignol and Amitsur showed that if $A$ is an
{\em abelian} splitting group of $\UD(p^r)$ and 
\begin{equation} \label{e.A_p}
A_p\isomo\bbZ/p^{n_1}\bbZ\times\dots\times\bbZ/p^{n_l}\bbZ 
\end{equation}
is its Sylow $p$-subgroup
then $n_\nu+n_{\nu+1}\ge[r/\nu]$ for every $\nu = 1, 2, \dots$;
see~\cite[Theorem 7.4]{ta1}. Consequently, the order of $A_p$ (and, hence,
of $A$) is divisible by $p^{f(r)}$, where
\[ f(r) = r+\sum_{\nu\ge3}\left\{\frac{[r/\nu]}{2}\right\}\ ;  \]
see~\cite[Theorem 7.5]{ta1}.
Here $[x]$ is the greatest integer $\leq x$ and $\{x\}$
is the smallest {\em nonnegative} integer $\geq x$. Note that 
$f(r) = \frac{1}{2}r\ln(r) + O(r)$, as $r \rightarrow \infty$;
see~\cite[Remark 7.5]{ta1} (a more precise asymptotic estimate is given 
in~\cite[Corollary 6.2]{ta2}). 

The same assertions hold if $A$ is an abelian splitting group
for any prime-to-$p$ extension of $\UD(p^r)$: the proof 
given in~\cite[Theorem 7.5]{ta1} goes through unchanged, except
that we use Proposition~\ref{prop.lagr} (with $e = 0$) in place of
\cite[Corollary 7.2]{ta1}, 

Moreover, let $D = \UD(p^r) \otimes_{Z(p^r)} K$, where $p^e$ is 
the highest power of $p$ which divides $[K:Z(p^r)]$, as in the statement 
of Theorem~\ref{thm1.3}. Suppose $A$ is a splitting group of $D$ and
$A_p$ is the Sylow $p$-subgroup of $A$. If $A_p$ is as in \eqref{e.A_p}
then a slight modification of the proof of~\cite[Lemma~6.1]{ta2} (again,
based on Proposition~\ref{prop.lagr}) shows that 
$n_\nu+n_{\nu+1}\ge[r/\nu] -e$ for every $\nu = 1, 2, \dots$
and consequently, the order of $A_p$ (and, hence,
of $A$) is divisible by $p^{f_e(r)}$, where
\[ f_e(r) =  r-e+\sum_{\nu\ge3}\left\{\frac{[r/\nu]-e}{2}\right\}\ . \]
It is easy to see that for a fixed $e$ and large $r$,
$f_e(r)$ also grows as $\frac{1}{2}r\ln(r) + O(r)$.
\end{remark}

\subsection*{Reduction of Theorem~\ref{thm.trdeg6} to a geometric problem}

Our proof of Theorem~\ref{thm.trdeg6} will be based 
on Proposition~\ref{prop6.8}. The idea is to construct a generically free 
$\PGL_{p^r}$-variety $X$ with two smooth points $x_1$ and $x_2$ 
whose stabilizers contain ``incompatible" symplectic modules $H_1$ and $H_2$.
Let $P_1$ and $P_2$ be as in \eqref{e.p1-p2}; this time we take
\[ H_2 = \phi_{P_2}(P_2 \times P_2^*) \simeq (\bbZ/p^r \bbZ)^2 \; , \]
as in \eqref{e.h1-h2}, but allow $H_1$ to be smaller:
\[ H_1 \,  = \, \text{rank 6 symplectic subgroup of} \; \phi_{P_1}(P_1 \times P_1^*) \;.\]
Note that $H_1 \simeq (\bbZ/p\bbZ)^6$ with desired properties exists 
by Lemma~\ref{lem.sub}(c).

Suppose $X$ is an irreducible generically free $\PGL_n$-variety,
and $x_1$, $x_2$ are smooth points of $X$ such that $x_i$ is fixed by 
$H_i$. Let $D$ be the algebra $RMaps_{\PGLn}(X, \Mn)$. Since $X$ has a smooth
point fixed by $H_2$, Lemma~\ref{lem.index} tells us that
$D$ is a division algebra.  Moreover, in view of Proposition~\ref{prop6.8} 
(with $X = X'$) any splitting group $A$ of $X$ (or equivalently, 
of $D$) will contain subgroups $L_1$ and $L_2$ which are isomorphic to
Lagrangian submodules
of $H_1$ and $H_2$, respectively.  Note that $L_1 \simeq (\bbZ/p\bbZ)^3$ and 
$L_2 \simeq (\bbZ/p^i \bbZ) \times (\bbZ/p^{r-i} \bbZ)$, for 
some $0 \leq i \leq r$. Then $L_1 \cap L_2$ is an abelian group of 
exponent $p$ and rank $\leq 2$. Thus 
\[ |A| \geq |L_1 L_2| = \frac{|L_1| \, |L_2|}{|L_1 \cap L_2|} \geq  
\frac{p^3 p^r}{p^2} = p^{r+1} \; ; \]
this shows that $D$ is not a a crossed product. The same argument
shows that any prime-to-$p$ extension of $D$ is not a crossed product.

Thus, in order to prove Theorem~\ref{thm.trdeg6} it is sufficient to
construct an irreducible generically free $\PGL_n$-variety $X$ such 
that $\trdeg_k k(X)^\PGLn = \dim(X/\PGL_n) = 6$
and $X$ has smooth points $x_1$
and $x_2$ such that $x_i$ is fixed by $H_i$.  

Note that both $H_1$ and $H_2$ are contained in the
finite subgroup $G$ of $\PGL_n$ generated by the permutation matrices 
and by the diagonal matrices all of whose entries are $p^r$th roots of unity.
We will construct $X$ as $\PGL_n \ast_{G} Y$, where $Y$ is
a 6-dimensional primitive $G$-variety with two points, 
$y_1$ and $y_2$ such that
$H_i$ fixes $y_i$. Indeed, if $Y$ is as above 
then the points $x_1 = (1_{\PGLn}, y_1)$ and 
$x_2 = (1_{\PGLn}, y_2)$ of $X$ have the desired properties. 
(Recall that $\PGLn \ast_G Y$ is defined as the geometric quotient of
$\PGLn \times Y$ by the $G$-action given by 
$g \cdot (h, y) = (hg^{-1}, g y)$; see~\cite[Section 4.8]{pv}.)

Therefore, in order to prove Theorem~\ref{thm.trdeg6} it is enough
to establish the following result.

\begin{prop} \label{prop.trdeg6} Let $G$ be a finite group and let
$H_1, \dots, H_s$ be abelian subgroups of $G$,
$r_i=\rank(H_i)$ and
$r = \max \{ r_i \mid i = 1, \dots, s \}$. Then there exists
a generically free primitive $r$-dimensional projective $G$-variety
$Y$ with smooth
points $y_1,\dots,y_s$ such that $H_i \subset \Stab(y_i)$.
\end{prop}

\begin{remark} \label{rem.jun22a}
Note that $\dim(Y)$ cannot be less than $r$.
More precisely, if $H$ is a finite abelian group, $Y$ is 
a quasiprojective $H$-variety, and $y$ is a smooth point of $Y$ fixed by $H$
then 
\begin{equation}\label{eqn.jun22a}
\codim_y(Y^{H})\ge\rank(H)\ .
\end{equation}
Indeed, assume the contrary: $\codim_y(Y^{H}) < \rank(H)$.  
By Lemma~\ref{lem.quasi-proj}(a), $y$ has an $H$-invariant affine
neighborhood in $Y$. Replacing $Y$ by this neighborhood,  
we may assume $Y$ is affine.
By the Luna Slice Theorem \cite[Corollary to Theorem~6.4]{pv}, $Y^H$
is smooth at $y$ and
\[
\dim(T_y(Y))-\dim(T_y(Y)^H)=\codim_y(Y^{H})< \rank(H)\ ;
\]
hence, the action
of $H$ on $T_{y}(Y)$ cannot be faithful. In other words, 
there exists a subgroup $H'\subset H$, $H'\ne\{1\}$, which acts trivially 
on $T_{y}(Y)$.
Applying \cite[Corollary to Theorem~6.4]{pv} to the action of $H'$ on
$Y$, we see that $H'$
acts trivially on all of $Y$. 
This contradicts our assumption that the $G$-action on $Y$ is generically 
free.  \qed
\end{remark}

\section{Constructing a $G$-variety with prescribed stabilizers}
\label{sect.trdeg6}

As we have just seen, Theorem~\ref{thm.trdeg6} follows from
Proposition~\ref{prop.trdeg6}.  This section 
will thus be devoted to proving Proposition~\ref{prop.trdeg6}.
Our general approach
is to first construct a higher-dimensional variety with desired
properties (this is easy), then replace it by a ``generic" $G$-invariant 
hypersurface passing through $y_1, \dots, y_s$, thus reducing 
the dimension by 1. To carry out this program, we first reduce to 
a situation where $Y^{H_i}$ has the highest possible dimension at $y_i$
(Theorem~\ref{thm.jun10a}), then apply Theorem~\ref{thm.bertini},
which may be viewed as a weak form of Bertini's theorem in the equivariant 
setting.

\subsection*{A local system of parameters}

The following lemma summarizes some known facts about
the local geometry of a smooth $G$-variety near
a point fixed by a finite abelian group. 

\begin{lem} \label{lem.jul23a}
Let $H$ be a finite abelian group, let $X$ be a smooth quasiprojective
$H$-variety,
and let $D_1,\dots,D_l$ be $H$-invariant hypersurfaces passing
through a point $x\in X$ and intersecting transversely at $x$.

\noindent
(1) There exists a local coordinate system (a regular system of
parameters) $u_1,\dots,u_n$ with the following properties:

(i) The group $H$ acts on each $u_i$ by a character $\xi_i$;

(ii) $u_i$ is the local equation of $D_i$ for $i=1,\dots,l$;

(iii) The germ of the fixed-point set $X^H$ at $x$ is given by the
local equations $u_{i_1}=\dots=u_{i_t}=0$ where $\{i_1,\dots,i_t\}$
is the set of all subscripts $i$ for which the character $\xi_i$ is
nontrivial.

\noindent
(2) Let $\pi\colon X'\to X$ be the blowup with the center $Z$ given by
the local equations $u_{j_1}=\dots=u_{j_s}=0$ at $x$.
(In particular, we can take $Z=X^H$.)
Let $\tilde D_i\subset X'$ be the strict transform of $D_i$.
Then we have:

(i) $\tilde D_1,\dots,\tilde D_l$ and the exceptional
divisor $\pi^{-1}(Z)$ are in normal crossing in a neighborhood of
$\pi^{-1}(x)$;

(ii) the natural isomorphism $\pi^{-1}(x) \isomo \bbP(T_xX/T_xZ)$
identifies $\tilde D_i\cap\pi^{-1}(x)$ with $\bbP(L_i)$, where $L_i$
is an $H$-invariant subspace of $T_xX/T_xZ$ of codimension 0 or 1.
\end{lem}

Note that if $u_i$ is a local equation of $D_i$ and $H$ acts on $u_i$
by a character $\xi_i$ then $H$ acts by the character $\xi_i$ on the
conormal space $(T_xX/T_xD_i)^*$.

\begin{pf} By Lemma~\ref{lem.quasi-proj}(a), we may assume without 
loss of generality that $X$ is affine.

(1) Denote by $\cO_x$ the local ring of $X$ at $x$, by $\fm_x$ its maximal
ideal, and  by $\fp_{D_i}$ the ideal of $D_i$ in $\cO_x$.

To construct $u_1,\dots,u_l$, note that the group $H$ acts on $\cO_x$,
the ideals $\fm_x$ and $\fp_{D_i}$ are $H$-invariant subspaces in $\cO_x$,
and the $H$-representation
$(\fp_{D_i}+\fm_x^2)/\fm_x^2\isomo\fp_{D_i}/(\fp_{D_i}\cap\fm_x^2)$ is
one-dimensional; let $\xi_i$ be the corresponding character of $H$.
The $H$-linear epimorphism
$\fp_{D_i}\to\fp_{D_i}/(\fp_{D_i}\cap\fm_x)$ splits; this yields a
generator $u_i\in\fp_{D_i}$ --- a local equation of $D_i$ --- on which
$H$ acts by the character $\xi_i$.

To construct $u_{l+1},\dots,u_n$, consider the $H$-linear epimorphism
\[
\fm_x\to\fm_x/(\fp_{D_1}+\dots+\fp_{D_l}+\fm_x^2)\isomo
\frac{\fm_x/\fm_x^2}{\sum_{i=1}^l(\fp_{D_i}+\fm_x^2)/\fm_x^2}\ ;
\]
its splitting yields the elements $u_{l+1},\dots,u_n\in\fm_x$ such
that $H$ acts on each of them by a character and the images of
$u_1,\dots,u_n$ in $\fm_x/\fm_x^2$ form a basis there.
It follows that $u_1,\dots,u_n$ form a regular system of parameters at
$x$ that satisfies properties (1)(i) and (1)(ii).

According to the Luna Slice Theorem \cite[Corollary to Theorem~6.4]{pv},
$X^H$ is given  in a neighborhood of $x$
by the local equations 
$u_{i_1}=\dots=u_{i_t}=0$, where $\{i_1,\dots,i_t\}$
is the set of all subscripts $i$ for which the character $\xi_i$ is
nontrivial. This proves part (1)(iii).

(2) Let $U$ be a small enough affine neighborhood of $x$ in $X$ so
that $u_1,\dots,u_n$ form a local coordinate system (i.e., their
differentials are linearly independent) everywhere on $U$.
The blown-up variety $X'$ in a neighborhood of $\pi^{-1}(x)$ is
covered by the charts $U_i$, $1\le i\le s$, where $U_i$ is
the complement in $\pi^{-1}(U)$ of the strict transform
of the subvariety $u_{j_i}=0$; the local coordinates in $U_i$ are
$v_{j_i}=\pi^*u_{j_i}$, $v_{j_{i'}}=\pi^*u_{j_{i'}}/\pi^*u_{j_i}$ for
$i'\ne i$, and
$v_j=\pi^*u_j$ for $j\not\in\{j_1,\dots,j_s\}$.  The exceptional
divisor in $U_i$ is given by the local equation $v_{j_i}=0$, and
$\tilde D_j$ is given by the local equation
$v_j=0$ in case $j\ne j_i$, and is empty  in case $j=j_i$.
Since the local equations of $\tilde D_1,\dots,\tilde D_l$
and of the exceptional divisor are elements of the same local
coordinate system, they are transverse; this proves (2)(i).

The local equations in $U_i$ of the preimage $\pi^{-1}(x)$ are $v_j=0$
for $j\not\in\{j_1,\dots,j_s\}$, and $v_{j_i}=0$; hence, $\pi^{-1}(x)$ is
contained in $\tilde D_j$ if $j\not\in\{j_1,\dots,j_s\}$, so that in
this case $\tilde D_j\cap\pi^{-1}(x)=\pi^{-1}(x)$ as claimed in (2)(ii).
If $j\in\{j_1,\dots,j_s\}$ then $D_i\supset Z$ and
$\tilde D_j\cap\pi^{-1}(x)$ can be identified with
$\bbP(T_xD_i/T_xZ)$; here $L_i=T_xD_i/T_xZ$ is an $H$-invariant subspace
of $T_xX/T_xZ$ of codimension 1 as claimed.
This completes the proof of (2)(ii).
\end{pf}

\begin{lem} \label{lem:jul27a} 
Let $G$ be an algebraic group, $H$ be a finite abelian subgroup of $G$, 
$X$ be a smooth quasiprojective $G$-variety and $x \in X^H$. Then there 
exists a sequence of blowups
\[ f \colon  X_i \lra \dots \lra X_0 = X \]
with smooth $G$-invariant centers, a point $y \in X_i^H$ satisfying
$f(y) = x$, and smooth $G$-invariant hypersurfaces $D_1,\dots,D_l$
meeting transversely, such that locally at $y$, the fixed point set
$X_i^H$ coincides with $D_1\cap\dots\cap D_l$.
\end{lem}

\begin{pf}
Let $\dim X=n$, and let
\[
X_n@>\pi_n>>\dots @>\pi_{i+1}>> X_i@>\pi_i>>\dots @>\pi_1>>X_0=X
\]
be the sequence of blowups with the centers $Z_i=X_i^H\subset X_i$.

If $X_i$ is smooth, then $Z_i=X_i^H$ is smooth by Luna's slice
theorem (see \cite[Corollary to Theorem~6.4]{pv}), and consequently,
$X_{i+1}$ is smooth. Thus $X_i$ and $Z_i$ are smooth for every $i$.
Note that the Luna Slice Theorem can be applied to the $G$-action on
$X_i$ because $X_i$ is quasiprojective and, hence, every point of $X_i$
has an open affine $G$-invariant neighborhood; 
see Lemma~\ref{lem.quasi-proj}.

Let $E_i = D_{i1}\cup\dots\cup D_{ii}$ be the exceptional divisor in $X_i$,
where $D_{ij}\subset X_i$ is the strict transform
(in $X_i$) of the exceptional divisor of $\pi_j\colon X_j\to X_{j-1}$.
We claim that $E_i$ is a normal crossing divisor. The proof is by induction
on $i$. The base case, $i = 0$, is obvious, since $E_0$ = the empty divisor, 
is normal crossing. For the inductive step, we assume that $E_i$
is a normal crossing divisor. Then $E_{i+1}$ is also a normal crossing divisor
by Lemma~\ref{lem.jul23a}(2)(i). This completes the proof of the claim. 

To obtain the required point $y \in X_m^H$, we 
start with $x_0 = x$ and inductively construct
$x_{i}\in X_i$ satisfying $\pi_{i}(x_{i})=x_{i-1}$
and $x_{i}\in X_i^H\cap D_{i1}\cap\dots\cap D_{ii}$, until 
we get a point $y = x_m$ with the desired properties.

Suppose $x_i$ has been constructed for some $i \geq 0$.

\smallskip
Case 1. The germ of $X_{i}^H$ at $x_{i}$ does not contain the
germ of $S=\bigcap_{j=1}^{i}D_{ij}$.
Let
\[ W=\pi_{i+1}^{-1}(x_{i})\cap\bigcap_{j=1}^{i+1}D_{i+1,j}\subset
X_{i+1}\ . \]
We claim $W \neq \emptyset$. Indeed, since $x_{i}\in X_{i}^H=Z_i$, we have
$\pi_{i+1}^{-1}(x_{i})\subset\pi_{i+1}^{-1}(Z_i)=D_{i+1,i+1}$, and thus
\begin{equation} \label{e.case1}
W=\pi_{i+1}^{-1}(x_{i})\cap\bigcap_{j=1}^{i}D_{i+1,j} \ .
\end{equation}
Since $D_{i+1, j}$ is the strict transform of $D_{ij}$,
$\bigcap_{j=1}^{i}D_{i+1,j}$ contains the strict transform of 
$S=\bigcap_{j=1}^{i}D_{ij}$. Thus $W$ contains the
intersection of the strict transform of $S$ with
$\pi_{i+1}^{-1}(x_{i})$.  As the germ of $S$ at
$x_{i}$ is not contained in the germ of the blowup center
$Z_i=X_{i}^H$, the strict transform of $S$ is nonempty and
intersects $\pi_{i+1}^{-1}(x_{i})$. Consequently, $W$ is nonempty, 
as claimed.

We now identify $\pi_{i+1}^{-1}(x_{i})$ with
$\bbP(T_{x_{i}}X_{i}/T_{x_{i}}Z_i)$ in the usual manner. 
Note that this identification is $H$-equivariant.
Then by Lemma~\ref{lem.jul23a}(2)(ii), 
$\pi_{i+1}^{-1}(x_{i})\cap D_{i+1,j}$ is identified with $\bbP(L_j)$,
where $L_j$ is an $H$-invariant subspace of the normal space 
$T_{x_{i}}X_{i}/T_{x_{i}}Z_i$.  Thus $W$ is identified with
$\bbP(L)$, where $L = L_1 \cap \dots \cap L_{i}$; see \eqref{e.case1}.  
Note that $L \neq (0)$ because $W \neq \emptyset$; moreover, $H$ acts
linearly on $L$. Since $H$ is diagonalizable, it has an eigenvector in $L$, 
i.e., a fixed point in $W \isomo \bbP(L)$. This fixed point is our new
point $x_{i+1}$.
By our construction, $\pi_{i+1}(x_{i+1})=x_i$ and
$x_{i+1}\in X_{i+1}^H\cap\bigcap_{j=1}^{i+1}D_{i+1,j}$, as desired.

\smallskip
Case 2: $X_{i}^H$ at $x_{i}$ contains $\bigcap_{j=1}^{i}D_{ij}$
at $x_{i}$.
Note that this is necessarily the case if $i=n$.
By Lemma~\ref{lem.jul23a}(1),
$X_{i}^H$ at $x_{i}$ coincides with the
intersection of those of $D_{ij}$ for which the action of $H$ on
$T_{x_{i}}X_{i}/T_{x_{i}}D_{ij}$ is nontrivial; in
particular, $X_{i}^H$ at $x_{i}$ is an intersection of smooth
$G$-invariant hypersurfaces meeting transversely, as required.

Thus, we see that for some $i\le n$, Case 2 occurs and we get a point
$x_i\in X_i$ with the required properties.
\end{pf}

\subsection*{Resolving the action on the tangent space}

In this subsection we prove the following result.

\begin{thm} \label{thm.jun10a}
Let $G$ be an algebraic group, let $X$ be a smooth quasiprojective $G$-variety,
let $H_1,\dots,H_s$ be finite abelian subgroups of $G$ and 
$x_1,\dots,x_s$ be points of $X$ such that $x_i$ is fixed by
$H_i$. Denote the rank of $H_i$ by $r_i$.
Then there 
is a sequence of blowups $\pi\colon X_m\to\dots\to X_0=X$ with smooth 
$G$-invariant centers, and points $y_1,\dots,y_s\in X_m$,
such that $\pi(y_i) = x_i$, $y_i\in X_m^{H_i}$, and
$X_m^{H_i}$ has codimension $r_i$ at $y_i$.
\end{thm}

Our proof relies on the following two simple lemmas.

\begin{lem} \label{lem.jun14a}
Let $H$ be a finite abelian group, let $X$ be an $H$-variety and
$x\in X^H$ be a smooth point of $H$.  If $\pi\colon\tilde X\to X$
is a blowup with a smooth $H$-invariant center $Z\subset X$
then there exists a point $\tilde x\in\pi^{-1}(x)$ fixed by $H$
such that $\dim_{\tilde x}(\tilde X^H)\ge \dim_x(X^H)$.
\end{lem}

\begin{pf}
Replacing $X$ by the set of its smooth points (which is clearly
$H$-invariant), we may assume that $X$ is smooth.  We claim that 
$\pi(\tilde{X}^H) = X^H$. The inclusion 
$\pi(\tilde{X}^H) \subset X^H$ is obvious. To prove the opposite 
inclusion, note that since $\pi$ is an isomorphism over
$X \setminus Z$, every $y \in X^H - Z$ lies in $\pi(\tilde{X}^H)$.
On the other hand, if $y \in Z^H$ then $\pi^{-1}(y)$ 
can be identified with $\bbP(T_yX/T_yZ)$ as $H$-varieties, and the
(linear) action of $H$ on $\bbP(T_yX/T_yZ)$ has a fixed point $\tilde{y}$; 
then $y = \pi(\tilde{y}) \in \pi(\tilde{X}^H)$.  This proves the claim, 
and the lemma follows.
\end{pf}

\begin{remark} Lemma~\ref{lem.jun14a} remains true under the more general
assumption that $H$ is Levi-commutative rather than finite abelian; see
Definition~\ref{def.levi} and Lemma~\ref{lem.levi}(iv).
The version we stated is sufficient for our application.
\end{remark}

\begin{lem} \label{lem:jul27b}
Let $A$ be a finite abelian group of rank $r$. An elementary operation
on $A^s$ is one of the form  \[ (\xi_1,\dots,\xi_i,\dots,\xi_s) \lra
(\xi_1,\dots,\xi_i-\xi_j,\dots,\xi_s) \]
for some $1\le i, j\le s$, where $i \neq j$.

Assume $s \ge r$. Then any
$\xi=(\xi_1,\dots,\xi_s) \in A^s$ can be transformed, by a finite sequence 
of elementary operations, into an $s$-tuple with at least $s-r$ zeros.
(Here by a zero, we mean the identity element of $A$.)
\end{lem}

\begin{pf} First we note that it does no harm to permute the components
of $\xi$. In other words, we may as well consider an operation of the form 
$(\xi_1, \dots, \xi_s) \mapsto (\xi_{\sigma(1)}, \dots, \xi_{\sigma(s)})$  
with $\sigma \in {\rm S}_n$, as another type of elementary operation.
The assertion of the lemma is then equivalent to saying that any $\xi \in A^s$ 
can be transformed, by these two types of elementary operations, into an
$s$-tuple $\lambda = (\lambda_1, \dots, \lambda_r, 0_A, \dots, 0_A)$, where
$0_A$ is the identity element of $A$.

We will prove this assertion by induction on $r$. Suppose $r = 1$, i.e., 
$A = \bbZ/n\bbZ$ for some $n \geq 1$. We can use elementary operations to
perform the Euclidean algorithm on $\xi_1$ and $\xi_2$. After 
interchanging them if necessary, we may assume $\xi_2 = 0$. (The
new value of $\xi_1$ is the greatest common divisor of the old 
values of $\xi_1$ and $\xi_2$). Applying the same procedure to $\xi_1$
and $\xi_3$, then $\xi_1$ and $\xi_4$, etc., we reduce the original
$s$-tuple to $(\xi_1, 0, \dots, 0)$, as claimed. 

For the induction step, write $A = B \times C$, where $B$ has
rank $r-1$ and $C$ is cyclic. Set
$\xi_i = (\beta_i, \gamma_i)$, where $\beta_i \in B$ and 
$\gamma_i \in C$. As we saw above, after performing a sequence of elementary
operations, we may assume $\gamma_2 = \dots = \gamma_s = 0_C$.
By the induction assumption, there exists a sequence of elementary operations
in $B^{s-1}$ which reduces $(\beta_2, \dots, \beta_s)$ to
$(\lambda_2, \dots, \lambda_r, 0_B, \dots, 0_B)$.
(Note that since $r \leq s$,
$\rank(B) = r-1 \leq s-1$, so that we may, indeed, use the induction
assumption.) Applying the same sequence to $(\xi_2, \dots, \xi_s)$, 
we reduce $(\xi_1, \dots, \xi_s)$ to 
\[ (\xi_1, (\lambda_2, 0_C), \dots, (\lambda_r, 0_C), 0_A, \dots, 0_A) 
\in A^s \; . \] 
This completes the proof of the lemma.
\end{pf}

\begin{pf*} {Proof of Theorem~\ref{thm.jun10a}}
By \eqref{eqn.jun22a}, we have
\begin{equation} \label{e9.1}
\codim_{x_i}(X^{H_i})\ge r_i \end{equation}
for any $i$.  We want to modify $X$ by a sequence of blowups so as 
to decrease $\codim_{x_i} X^{H_i}$ to $r_i$ for each $i$. (Of course, after
each blowup $\tilde{X} \lra X$ we replace $X$ by $\tilde{X}$ and $x_i$
by $\tilde{x_i}$, as in Lemma~\ref{lem.jun14a}.) We claim that
we may do this for one $i$ at a time; in other words, we may assume $s = 1$.  
Indeed, suppose we have reduced to the case where 
$\codim_{x_1}(X^{H_1}) = r_1$. If we now perform a further blowup
$\tilde{X} \lra X$ and choose $\tilde{x_1}$ above $x_1$ as in 
Lemma~\ref{lem.jun14a}, then  Lemma~\ref{lem.jun14a} and \eqref{e9.1} 
tell us that $\codim_{\tilde{x_1}}(\tilde{X}^{H_1}) = r_1$. Thus
we are free to perform another sequence of blowups that would give us the
desired equality for $i = 2$, then $i = 3$, etc.

We will thus assume $s=1$ and set $x = x_1$, $H = H_1$, $r = r_1$.

By performing a sequence of blowups given by Lemma~\ref{lem:jul27a},
we may assume that there exist 
$G$-invariant divisors $D_1, \dots, D_c$ such that
$X^H = D_{1}\cap\dots\cap D_{c}$ in a neighborhood of $x$, 
$c=\codim_{x}X^H$ and $D_1, \dots, D_c$ intersect at $x$
transversely. 

Note that $T_x(X)/T_x(X_i^H) \isomo \bigoplus_{j=1}^cT_x(X)/T_x(D_j)$. 
Here $H$ acts on each one-dimensional space
$T_x(X)/T_x(D_i)$ by a character $\xi_i\in H^*$ which is nontrivial by
Lemma~\ref{lem.jul23a}(1)(iii). In other words, the linear action
of $H$ on the tangent space $T_x(X)$ decomposes as a direct 
sum of $c$ nontrivial characters
$\xi_{1},\dots,\xi_{c}$ and $n-c = \dim X^H$ trivial characters.

Recall that by \eqref{e9.1}, $c\ge r$. We would like to modify $X$ by a
sequence of blowups to arrive at the situation where $c = r$. In other
words, if $c > r$, we want to perform a sequence of blowups that
would lower the value of $c$.

With this goal in mind, we would like to know how the characters $\xi_i$
change after one blowup. Specifically, we will consider the blowup
$\pi \colon \tilde{X} \lra X$ with center $Z = D_i \cap D_j$, where 
$1 \leq i,j \leq c$, $i\ne j$.  Since $Z$ is of codimension 2 in $X$,
$\pi^{-1}(x)$ is isomorphic to $\bbP^1$.
Let $\tilde{x}$ be the (unique) point of $\pi^{-1}(x)$
that lies in the strict
transform of $D_{i}$, take $\tilde{D}_l$ to be the strict transform of
$D_l$ for $l=1,\dots,\hat{j},\dots,c$, and let
$\tilde{D}_j=\pi^{-1}(Z)$ be the exceptional divisor of $\pi$.
Then the action of $H$ in $T_{\tilde{x}}\tilde X$ is given by the direct
sum of the characters
$\tilde{\xi_l}=\xi_{l}$ if $l\ne i$,
$\tilde{\xi_i}=\xi_{i}\xi_j^{-1}$,
and $(\dim X-c)$ trivial characters. In other words, the new characters
$\tilde{\xi_1},\dots, \tilde{\xi_c} \in H^{\ast}$ are obtained from
the old characters $\xi_{1},\dots,\xi_{c} \in H^{\ast}$ 
by an elementary operation, as in Lemma~\ref{lem:jul27b}.
(Note that our group $H^{\ast}$ is written multiplicatively, whereas
the group $A$ in Lemma~\ref{lem:jul27b} is written additively.)

Now Lemma~\ref{lem:jul27b} tells us that there is a sequence of
of elementary operations which transforms $(\xi_1, \dots, \xi_c)$
to $(\lambda_1, \dots, \lambda_r, 1_{H^{\ast}}, \dots, 1_{H^{\ast}})$
for some $\lambda_1, \dots, \lambda_r \in H^{\ast}$. Recall that initially
our characters $\xi_i$ are all nontrivial. We want to follow the above sequence
of elementary operations until we create the first trivial character.
Each one of these transformations is given by a blowup of a codimension 2
subvariety, as described above. When the first trivial character appears, 
$\dim(X^H)$ goes up by one. At that point, we reduce $c$ by 1 and
repeat the above procedure, until $c$ becomes equal to $r$.
\end{pf*}

\subsection*{A weak equivariant Bertini theorem}

\begin{thm} \label{thm.bertini}
Let $G$ be a finite group, $X$ a primitive smooth projective $G$-variety,
$x_1,\dots,x_s$ points of $X$ with stabilizers
$H_i=\Stab(x_i)$ for $i=1,\dots,s$.
Suppose that each $x_i$ is not an isolated point of $X^{H_i}$ and
$\dim(X)\ge2$.
Then 

\smallskip
(a) there exists a smooth closed $G$-invariant primitive hypersurface
$W\subset X$ passing through $x_1,\dots,x_s$.

\smallskip
(b) Moreover, if $X$ is a generically free $G$-variety then we can
choose $W$ so that it is also generically free.  
\end{thm}

We shall need the following variant of Bertini's theorem; for lack
of a reference we will supply a proof.

\begin{lem} \label{lem:gen}
Let $x_1,\dots,x_s$ be points in $\bbP^n$ and let
$V_d$ be the space of homogeneous polynomials of degree $d\ge s+1$ in
$\bbP^n$ that vanish at $x_1,\dots,x_s$.

(a) Suppose $Y$ is a locally closed subvariety of $\bbP^n$, 
$y$ is a smooth point of $Y$ and
{\normalshape $ 
V_{d,y}=
\Bigl\{P\in V_d\Bigm| P|_Y\text{ has zero of order }>1\text{ at
}y\Bigr\}$}.
Then the codimension of $V_{d,y}$ in $V_d$ is given by
{\normalshape\[ \codim(V_{d,y})=
\begin{cases}
\dim_y(Y)+1 &\text{ if }y\not\in\{x_1,\dots,x_s\}\ ,\\
\dim_y(Y) &\text{ if }y\in\{x_1,\dots,x_s\}\ .
\end{cases} \]}

(b) Suppose $Y_1,\dots,Y_l$ are smooth locally closed subvarieties of
$\bbP^n$ such that 
\begin{equation} \label{e.lem9.8}
 \text{{\normalshape $x_j$ is not an isolated point of $Y_i$ for all $i,j$.}}
\end{equation}
Then a generic hypersurface of degree $d\ge s+1$ in $\bbP^n$ that passes
through $x_1,\dots,x_s$, is transverse to $Y_1,\dots,Y_l$.
\end{lem}

Note that assumption \eqref{e.lem9.8} in part (b) is necessary, since 
otherwise no hypersurface passing through $x_j$ can be transverse to $Y_i$.  
(By definition, a hypersurface $W$ is transverse to a one-point 
set $\{x_j\}$ iff $W$ does not pass through $x_j$.)

\begin{pf} (a) Choose an affine subset
$\bbA^n=\Spec k[z_1,\dots,z_n]\subset\bbP^n$ that contains
$x_1,\dots,x_s$ and $y$; then $V_d$ may be identified with the space
of all polynomials in $z = (z_1,\dots,z_n)$ of degree $\le d$ (not necessarily
homogeneous) that vanish at the points $x_1,\dots,x_s$.

Consider the linear map $\phi_y\colon V_d\to\cO_y/\fm_y^2$, where
$\cO_y/\fm_y^2$ is the space of $1$-jets of regular functions on
$\bbA^n$ at $y$;
$\phi_y$ sends each polynomial $P\in V_d$ into its $1$-jet at $y$
(cf.\ \cite[Proof of Theorem~8.18]{Hart}).
Then $V_{d,y}=\phi_y^{-1}(N_y)$ where $N_y$ is the subspace of
$\cO_y/\fm_y^2$ consisting of the jets of all functions that vanish on
$Y$; $N_y$ may be identified with the conormal space to
$Y$ at $y$, and hence, it is a linear subspace of
$\cO_y/\fm_y^2$ of dimension $\dim(N_y)=\codim_y(Y)$.

Assume $y\not\in\{x_1,\dots,x_s\}$. 
Let $l(z)$ be a linear combination
of $z_1,\dots,z_n$ whose value at $y$ is different from its values at
$x_1,\dots,x_s$.  
Then the degree $s$ polynomial
$P(z) = \Bigl(l(z) - l(x_1)\Bigr) \dots \Bigl(l(z) - l(x_s)\Bigr)$ 
vanishes at
$x_1,\dots,x_s$ and does not vanish at $y$. This means that
$\phi_y(P(z))$ and $\phi_y(z_jP(z))$ (where $j=1,\dots,n$)
span  $\cO_y/\fm_y^2$ as a $k$-vector space.
Hence, $\phi_y$ is onto, and
\begin{multline*}
\codim(V_{d,y})=\codim(\phi_y^{-1}(N_y))=\dim(\cO_y/\fm_y^2)-\dim(N_y)\\
=n+1-\codim_y(Y)=\dim_y(Y)+1
\end{multline*}
as claimed.

Now suppose $y\in\{x_1,\dots,x_s\}$, say, $y=x_1$.
In this case
$\phi_y(V_d)$ is clearly contained in $\fm_y/\fm_y^2$; we will show
that, in fact, equality holds. Indeed, we may assume without loss of
generality that $y = x_1 = (0, \dots, 0) \in \bbA^n$. Let $l(z)$ be a
linear combination of $z_1, \dots, z_n$ such that $l(x_2), \dots, l(x_s) \neq
l(x_1) = 0$ and let
$Q(z) = \Bigl(l(z)-l(x_2)\Bigr)\dots \Bigl(l(z) - l(x_s)\Bigr)$.
Given any linear function $a(z)$ in $z_1, \dots, z_n$, let 
$R(z) = a(z)Q(z)$. Note that $R(z)$ vanishes at $x_1, \dots, x_s$ 
and has degree $s$; hence, $R(z) \in V_d$. On the other hand,
$\phi_y(R)$ equals $R(z)$ modulo the terms of degree $\geq 2$ in $z_1, \dots,
z_n$, i.e., $\phi_y(R) = Q(0)a(z)$ (mod $\fm_y^2$), 
where $Q(0)$ is a nonzero element of
$k$. This means that $\phi_y(V_d)$ contains $a(z)$, thus proving that
$\phi_y(V_d) = \fm_y/\fm_y^2$, as desired. Now
\begin{multline*}
\codim(V_{d,y})=\codim(\phi_y^{-1}(N_y))=\dim(\fm_y/\fm_y^2)-\dim(N_y)\\
=n-\codim_y(Y)=\dim_y(Y)\ .
\end{multline*}

(b) After replacing each $Y_i$ by the collection of its irreducible components, 
we may assume each $Y_i$ is irreducible. Moreover, we may assume without
loss of generality that $l = 1$; we shall denote $Y_1$ by $Y$.

Let $X$ be the algebraic subvariety of $V_d \times Y$ given by
\[ X = \Bigl\{(P, y) \Bigm| P|_Y\text{has a zero of order }>1\text{ at } 
y\Bigr\}\ ; \]
Denote the natural projections of $X$ to $V_d$ and $Y$ by $\pi_1$
and $\pi_2$. We want to show that $\dim(\pi_1(X)) < \dim(V_d)$.
It is thus enough to prove that $\dim(X) < \dim(V_d)$. Write 
$X = X_1 \cup X_2$, where
$X_1 = \pi_2^{-1}(Y - \{x_1, \dots, x_s \})$ and  
$X_2 = \pi_2^{-1}(Y \cap \{x_1, \dots, x_s \})$. It is enough
to show that $\dim(X_i) < \dim(V_d)$ for $i = 1, 2$. 
The fibers of $\pi_2$ are precisely the sets $V_{d, y}$ we considered in
part (a). Since $\dim(V_{d, y}) = \dim(V_d) - \dim(Y) - 1$ for every
$y \in Y - \{x_1, \dots, x_s \}$, we conclude that
$\dim(X_1) \leq \dim(V_d) -1$.

It remains to show that $\dim(X_2) \leq \dim(V_d) -1$. If
$Y \cap \{x_1, \dots, x_s \} = \emptyset$ then $X_2 = \emptyset$, and
there is nothing to prove. On the other hand, if 
$Y \cap \{x_1, \dots, x_s \} \neq \emptyset$ then assumption \eqref{e.lem9.8}
says that $\dim(Y) \geq 1$.  Thus by part (a), 
\[ \dim(X_2) \leq \dim(V_d) - \dim(Y) \leq \dim(V_d) - 1 \;,  \]
as claimed.
\end{pf}

\begin{pf*}{Proof of Theorem~\protect\ref{thm.bertini}}
We begin with three simple observations. First of all,
we may assume without loss of generality that
the orbits $Gx_i$ are disjoint. Indeed, if $W$ passes through 
$x_i$ then it will pass through every point of $Gx_i$. Thus if, say,
$x_j$ happens to lie in $Gx_i$ then we can simply remove $x_j$ from 
our finite collection of points and proceed to construct $W$ for the
smaller collection. 

Secondly, part (b) is an immediate consequence of part (a). Indeed, since
$G$ is a finite group, generically free $G$-varieties are precisely
faithful $G$-varieties, i.e., $G$-varieties, where every nonidentity element
of $G$ acts nontrivially. The set
\[
X_0 = \{ x \in X \mid \Stab(x) = \{ 1 \} \}
\]
is open and dense in $X$;
in order to ensure that $W$ is generically free, it is enough to 
construct $W$ so that $W \cap X_0 \neq \emptyset$. This is accomplished by
applying part (a) to the collection $\{ x_0 , x_1, \dots, x_s \}$, where
$x_0 \in X_0$. Therefore, it is enough to prove part (a).

Thirdly, since $G$ is a finite group and $X$ is projective, there exists
a (finite) geometric quotient morphism $\psi\colon X\to X//G$
with $X//G$ projective; see Lemma~\ref{lem.quasi-proj}.
Since $X$ is a primitive $G$-variety, $X//G$ is irreducible. (Recall 
that the geometric quotient $X//G$ is a birational model for the rational
quotient $X/G$ which is irreducible since $X$ is primitive.) Note that
$X$
is partitioned into a union of nonintersecting smooth locally
closed subsets $\tilde{X^H} = \{ x \in X \, | \, \Stab(x) = H \}$, where
$H$ ranges over the set of subgroups of $G$.
By the Luna Slice Theorem~\cite[Theorem~6.1]{pv} the morphism
\begin{equation} \label{eqn.jun22b}
\psi|_{\tilde{X^H}}\colon\tilde{X^H}\to\psi(\tilde{X^H})
\end{equation}
is \'etale, and hence, the sets $\psi(\tilde{X^H})$ are also smooth.
(Note that Luna's theorem can be applied to
the $G$-action in a neighborhood of any point of $X$ by
Lemma~\ref{lem.quasi-proj}(a).) 
Two subvarieties $\psi(\tilde{X^H})$ and $\psi(\tilde{X^{H'}})$
coincide if the subgroups $H$ and $H'$ are conjugate, and are disjoint
otherwise.
In other words, $X//G$ is partitioned into a union of nonintersecting
smooth locally closed subsets $\psi(\tilde{X^H})$ for all conjugacy 
classes of subgroups $H\subset G$, so that every point $y\in X//G$
lies in $\psi(\tilde{X^H})$ where $H$ belongs to the conjugacy 
class of the stabilizers of points in $\psi^{-1}(y)$.

We are now ready to proceed with the proof of part (a).
Since $X//G$ is projective, we can embed it in $\bbP^N$ for some
$N$. Let
$U$ be a generically chosen hypersurface
of degree $s+1$ in $\bbP^N$ which passes through
$\psi(x_1),\dots,\psi(x_s)$, and let
$\overline{W} = X//G \cap U$.  We claim that 
\[ W = \pi^{-1}(\overline{W}) \subset X \]
satisfies the conditions of part (a). By our construction $W$ is 
a closed $G$-invariant hypersurface in $X$ passing through $x_1, \dots, x_s$;
thus we only need  to show that $W$ is smooth and primitive.

Since each $x_i$ lies in $X^{H_i}$ and is not its isolated point,
each $\psi(x_i)$ lies in $\psi(\tilde{X^{H_i}})$ and is not its
isolated point; it follows that $\psi(x_i)$ is not an isolated point
of any $\psi(\tilde{X^{H}})$.
By Lemma~\ref{lem:gen}(b), $U$ intersects every
subvariety $\psi(\tilde{X^H})$ transversely.

Let $x\in W$ and $H=\Stab(x)$; then $x\in\tilde{X^H}$ and
$\psi(x)\in\psi(\tilde{X^H})$.
Let $f$ be a local equation of $U$ at $\psi(x)$; then $\psi^*(f)$ is a
local equation of $W$ in $X$ at $x$.
Since $U$ is transverse to $\psi(\tilde{X^H})$, the restriction
$f|_{\psi(\tilde{X^H})}$ is nondegenerate at $\psi(x)$.
As the morphism \eqref{eqn.jun22b} is \'etale,
$\psi^*(f|_{\psi(\tilde{X^H})})=\psi^*(f)|_{\tilde{X^H}}$ is a
nondegenerate function on $\tilde{X^H}$ at $x$.
Hence, $\psi^*(f)$ is nondegenerate at $x$; in other words,
$W$ is smooth at $x$.

It remains to prove that $W$ is primitive or, equivalently, $\overline{W}$
is irreducible.  By~\cite[Corollary~III.7.9]{Hart}, 
$\overline W = X//G \cap U$ is connected.
On the other hand, since $W$ is smooth, $\overline W=W//G$ is normal.
We conclude that $\overline W$ is irreducible, as claimed.
\end{pf*}

\subsection*{Proof of Proposition~\ref{prop.trdeg6}}
If $r=0$ then $H_1=\dots=H_s=\{1\}$, and we can take
$Y$ to be a set of $|G|$ points with a transitive $G$-action
and $y_1, \dots, y_s$ to be any $s$ points in $Y$ (not necessarily 
distinct).  From now on we shall assume $r\ge1$.

Let $V$ be a generically free linear representation of $G$. By
Proposition~\ref{prop.proj} and Theorem~\ref{thm.resol}, there exists
a smooth projective $G$-variety $X$ such that $X \simeq V$ (as $G$-varieties) 
and $\Stab(x)$ is commutative for any $x\in X$. 
(Note that a Levi-commutative finite group is
commutative.) Every $H_i$ has a smooth fixed point in $V$, namely
the origin. Applying 
the Going Down Theorem~\ref{going-down} to the birational isomorphism
$V \stackrel{\simeq}{\brokrarr} X$, we conclude that $X^{H_i} \neq \emptyset$
for every $i$. The resulting smooth projective irreducible generically free
$G$-variety $X$ is the starting point for our construction.

After birationally modifying $X$ by a sequence of blowups 
with smooth $G$-equivariant centers, we may assume that there are
points $x_1,\dots,x_s\in X$ such that each $x_i$ is fixed by $H_i$
and the codimension of $X^{H_i}$ at $x_i$ is $r_i$; see
Theorem~\ref{thm.jun10a}. 

If $\dim X>r=\max_ir_i$ then $\dim X>r_i=\codim_{x_i}X^{H_i}$ and
hence, $x_i$ is not an isolated fixed point of $H_i$ for each $i$.
In addition, $\dim X>r\ge1$ implies $\dim X\ge2$.  Then
Theorem~\ref{thm.bertini}(b) yields a smooth closed generically free
$G$-invariant
primitive hypersurface $W$ in $X$ passing through $x_1,\dots,x_s$.
Replacing $X$ by this hypersurface reduces $\dim X$ by one.
Applying this procedure $\dim X-r$ times, we obtain a smooth
$G$-invariant primitive
subvariety $Y$ of dimension $r$ passing through $x_1,\dots,x_s$,
and hence, having points fixed by $H_1,\dots,H_s$. This completes
the proof of Proposition~\ref{prop.trdeg6} and thus of 
Theorem~\ref{thm.trdeg6}.
\qed

\begin{remark} \label{rem.Q} A closer examination of the proof of 
Proposition~\ref{prop.trdeg6} shows that the $G$-variety 
$Y$ can, in fact, be constructed over $\bbQ$.
Thus the division algebra $D$ is Theorem~\ref{thm.trdeg6} can be assumed
to be defined over $\bbQ$. This means that there exists a finitely generated 
field extension $F/\bbQ$ and a division algebra $D_0$ with center $F$ such 
that $\trdeg_{\bbQ} (F) = 6$ and $D = D_0 \otimes_F K$.
\end{remark}

\begin{remark} \label{rem:thm1.5+}
Our argument can be modified to prove the following extension of
Theorem~\ref{thm.trdeg6}: for any integer $e\ge0$
there exists a division algebra $D$ with center $K$ such that

\nobreak
\smallskip
(a) $K$ is a finitely generated extension of $k$ of transcendence
degree $6+2e$ and

\smallskip
(b) any extension of $D$ of degree $s$ is not a crossed product, provided
that $p^{e+1}\notdiv s$.
\end{remark}

\end{document}